\documentclass[12pt]{article}
\usepackage{amsmath,amsthm,amsfonts,amssymb}
\def\qed{{\unskip\nobreak\hfil\penalty50
\hskip2em\hbox{}\nobreak\hfil$\square$
\parfillskip=0pt \finalhyphendemerits=0\par}\medskip}
\def\proof{\trivlist \item[\hskip \labelsep{\bf Proof\ }]}
\def\endproof{\null\hfill\qed\endtrivlist}

\def\isom{\cong}

\def\lan{\langle}
\def\ran{\rangle}
\def\Ad{{\hbox{\rm Ad}}}

\def\dim{{\hbox{dim}}}

\def\id{{\rm id}}
\def\o{{\rm opp}}

\def\a{\alpha}
\def\be{\beta}
\def\de{\delta}
\def\e{\varepsilon}
\def\epsilon{\varepsilon}
\def\ga{\gamma}

\def\phi{\varphi}
\def\si{\sigma}

\def\s{{\sigma}}
\def\p{{\pi}}

\newtheorem{theorem}{Theorem}
\newtheorem{lemma}[theorem]{Lemma}

\newtheorem{corollary}[theorem]{Corollary}
\newtheorem{definition}[theorem]{Definition}

\newtheorem{proposition}[theorem]{Proposition}

\def\emptyset{\varnothing}
\def\setminus{\smallsetminus}

\def\2#1{{\cal #1}}
\def\ov{\overline}
\def\A{{\cal A}}
\def\B{{\cal B}}
\def\C{{\cal C}}
\def\D{{\cal D}}
\def\O{{\cal O}}
\def\J{{\cal J}}
\def\I{{\cal I}}

\def\cN{{\cal N}}
\def\N{{\cal N}}
\def\M{{\cal M}}
\def\E{{\cal E}}
\def\H{{\cal H}}

\def\Z{{\mathbb Z}}
\def\C{{\mathbb C}}
\def\R{{\mathbb R}}

\def\ext{\,{\rm ext}\,}
\def\In{{\bf{I}}}
\def\AA{{\A\otimes \A^{\rm opp}}}

\title{{\bf
Multi-interval Subfactors \\
and Modularity of Representations\\
in Conformal Field Theory}}

\author{
{\sc Yasuyuki Kawahigashi}\\
Department of Mathematical Sciences\\
University of Tokyo, Komaba, Tokyo, 153-8914, Japan\\
e-mail: {\tt yasuyuki@ms.u-tokyo.ac.jp}\\
\vphantom{X}\\
{\sc Roberto Longo}\footnote{Supported in part by GNAFA and
MURST.} \phantom{i}and
{\sc Michael M\"uger}\footnote{Supported by EU TMR Network
``Noncommutative Geometry''.}\\
Dipartimento di Matematica\\
Universit\`a di Roma ``Tor Vergata''\\
Via della Ricerca Scientifica, I-00133 Roma, Italy\\
e-mail: {\tt longo@mat.uniroma2.it, mueger@x4u.desy.de}}
\begin{document}
\date{}
\maketitle
\centerline{\sl Dedicated to John E. Roberts on the occasion of his
sixtieth birthday}
\bigskip

\begin{abstract}
We describe the structure of the inclusions of factors $\A(E)\subset
\A(E')'$
associated with multi-intervals $E\subset\R$ for a local
irreducible net $\A$ of von Neumann algebras on the real line satisfying
the split property and Haag duality.
In particular, if the net is conformal and the subfactor has finite
index,
the inclusion associated with
two separated intervals is isomorphic to
the Longo-Rehren inclusion, which provides a
quantum double construction of the tensor category of superselection
sectors of $\A$. As a consequence,  the index of  $\A(E)\subset
\A(E')'$
coincides with the global index associated with all irreducible
sectors, the braiding symmetry associated with all sectors is
non-degenerate, namely the representations of $\A$
form a modular tensor category, and every
sector is a direct sum of sectors with finite dimension.
The superselection structure is generated by local data.
The same results hold true if conformal invariance is
replaced by strong additivity
and there exists a modular PCT symmetry.
\end{abstract}
\newpage

\section{Introduction}

This paper provides the solution to a natural problem in (rational)
conformal quantum field theory, the description of the structure of the
inclusion of factors associated to two or more separated intervals.

This problem has been considered in the past years, seemingly with different
motivations. The most detailed study of this inclusion so far has been done by Xu
\cite{X2} for the models given by loop group construction for $SU(n)_k$ \cite{W}.
In this case Xu has computed the index and the dual principal graph of the inclusions. 
A suggestion to study this inclusion has been made also in \cite[Section 3]{S}.
Our analysis is model independent, and will display new structures and a deeper
understanding also in these and other models.

Let $\A$ be a local irreducible conformal net of von Neumann algebras on $\R$, i.e. an
inclusion preserving map 
$$
I\mapsto\A(I)
$$
from the (connected) open intervals of $\R$ to von
Neumann algebras $\A(I)$ on a fixed Hilbert space.
One may define $\A(E)$ for an arbitrary set $E\subset\R$ as the
von Neumann algebra generated by all the $\A(I)$'s as $I$ varies in
the intervals contained in $E$. By locality $\A(E)$ and $\A(E')$
commute, where $E'$ denotes the interior of
$\R\setminus E$, and thus one obtains an inclusion
\[
\A(E)\subset\hat\A(E),
\]
where $\hat\A(E)\equiv\A(E')'$. If Haag duality holds,
as we shall assume\footnote{As shown in \cite{GLW},
one may always extend $\A$ to the
dual net $\A^d$, which is conformal and satisfies Haag duality.},
this inclusion is trivial if $E$ is an interval, but it is in general
non-trivial
for a disconnected
region $E$. We will explain its structure if
$E$ is the union of $n$ separated intervals, a situation that can be reduced
to the case $n=2$, namely $E=I_1\cup I_2$, where $I_1$ and $I_2$ are
intervals with disjoint closure, as we set for the rest of this
introduction.

One can easily realize that the inclusion $\A(E)\subset\hat\A(E)$
is related to the superselection structure of $\A$, i.e. to the
representation theory of $\A$, as charge transporters between
endomorphisms localized in $I_1$ and $I_2$ naturally live in
$\hat\A(E)$, but not in $\A(E)$.

Assuming the index $[\hat\A(E):\A(E)]<\infty$ and the split
property\footnote{This general property is satisfied, in particular,
if Tr$(e^{-\beta L_0})<\infty$ for all $\beta>0$, where $L_0$ is the
conformal Hamiltonian, cf. \cite{BDL,DRL}.},
namely that $\A(I_1)\vee\A(I_2)$ is naturally isomorphic to
$\A(I_1)\otimes\A(I_2)$, we shall show that indeed $\A(E)\subset\hat\A(E)$
contains all the information on the superselection rules.

We shall prove that in this case $\A$ is rational, namely there exist
only finitely many different irreducible sectors $\{[\rho_i]\}$ with
finite dimension and that $\A(E)\subset\hat\A(E)$ is isomorphic to the
inclusion considered in \cite{LR} (we refer to this as the LR
inclusion, cf. Appendix \ref{A1}), which is canonically associated with
$\A(I_1)$, $\{[\rho_i]\}$ (with the identification
$\A(I_2)\simeq\A(I_1)^\o$).
In particular,
\[
[\hat\A(E):\A(E)]=\sum_i d(\rho_i)^2,
\]
the global index of the superselection sectors.
In fact $\A$ will turn out to be rational in an even stronger sense,
namely there exist no sectors with infinite dimension, except the
ones that are trivially constructed as
direct sums of finite-dimensional sectors.

Moreover, we shall exhibit an explicit way to generate the superselection
sectors of $\A$ from the local data in $E$: we consider the canonical
endomorphism $\gamma_E$ of $\hat \A(E)$ into $\A(E)$ and its
restriction $\lambda_E = \gamma_E |_{A(E)}$; then $\lambda_E$ extends
to a localized endomorphism $\lambda$ of $\A$ acting identically on
$\A(I)$ for all intervals $I$ disjoint from $E$. We have
\begin{equation}
\lambda = \bigoplus_i \rho_i \bar\rho_i,
        \label{canonical}
\end{equation}
where the $\rho_i$'s are inequivalent irreducible endomorphisms of
$\A$ localized in $I_1$ with conjugates $\bar\rho_i$ localized in
$I_2$ and the classes $\{[\rho_i]\}_i$ exhaust all the irreducible
sectors.

To understand this structure, consider the symmetric case $I_1=I$, $I_2= -I$.
Then $\A(-I) = j(\A(I))$, where $j$ is the anti-linear PCT
automorphism, hence we may identify $\A(-I)$ with $\A(I)^\o$. Moreover
the formula $\bar\rho_i = j\cdot\rho_i\cdot j$ holds for the conjugate
sector \cite{GL}, thus by the split property we may identify
$\{\A(E), \rho_i \bar\rho_i |_{\A(E)}\}$ with
$\{\A(I)\otimes\A(I)^\o, \rho_i \otimes\rho_i^{\o}\}$.
Now there is an isometry $V_i$ that intertwines the identity and
$\rho_i\bar\rho_i$ and belongs to $\hat\A(E)$. We then have to show
that $\hat\A(E)$ is generated by $\A(E)$ and the $V_i$'s and that
the $V_i$'s satisfies the (crossed product) relations characteristic of
the LR inclusion. This last point is verified by  identifying $V_i$
with the standard implementation isometry as in \cite{GL}, while
the generating property follows by the index computation that will
follow by the ``transportability'' of the canonical endomorphism above.

The superselection structure of $\A$ can then be recovered
by formula (\ref{canonical}) and the split property.
Note that the representation tensor category of $\AA$ generated
by $\{\rho_i\otimes\rho_i^\o\}_i$ corresponds to the connected
component of the identity in the fusion graph for $\A$, therefore the
associated fusion rules and quantum $6j$-symbols
are encoded in the isomorphism class of
the inclusion $\A(E)\subset\hat \A(E)$, that will be completely
determined by a crossed product construction.

A further important consequence is that the braiding symmetry
associated with all sectors is always non-degenerate, in other words
the localizable representations form a modular tensor category. As shown by
Rehren \cite{R1}, this implies the existence and
non-degeneracy of Verlinde's matrices $S$ and $T$,
thus the existence of a unitary representation of the
modular group $SL(2,\mathbb Z)$, which plays a role in
topological quantum field theory.

It follows that the net $\B\supset\AA$ obtained by the LR construction
is a field algebra for $\AA$, namely $\B$ has no superselection sector
(localizable in a bounded interval) and there is a generating family
of sectors of $\AA$ that are implemented by  isometries in $\B$.
Indeed $\B$ is a the crossed product of $\AA$ by the tensor category of
all its sectors.

As shown by Masuda \cite{Ma},
Ocneanu's asymptotic inclusion \cite{O2}  and
the Longo-Rehren
inclusion in \cite{LR} are, from the categorical viewpoint,
essentially the same constructions. The construction of the asymptotic
inclusion gives a new subfactor
$\M\vee (\M'\cap \M_\infty)\subset \M_\infty$ from a hyperfinite
II$_1$ subfactor $\N\subset \M$ with finite index and finite
depth and it is a subfactor analogue of the quantum
double construction of Drinfel$'$d \cite{D}, as noted
by Ocneanu.  That is, the tensor category of the $\M_\infty$-$\M_\infty$
bimodules arising from the new subfactor is regarded a ``quantum
double'' of the original category of $\M$-$\M$ (or $\N$-$\N$) bimodules.

On the other hand, as shown in \cite{Mu3}, the Longo-Rehren
construction gives the quantum double of the original tensor category
of endomorphisms.
(See also \cite[Chapter 12]{EK2}
for a general theory of asymptotic inclusions and their
relations to topological quantum field theory.)

Our result thus shows that the inclusion
arising from two separated intervals as
above gives the quantum double of the tensor category of all
localized endomorphisms.
However, as the braiding symmetry is non-degenerate, the quantum
double will be isomorphic to the subcategory of the trivial doubling of
the original tensor category corresponding to the connected
component of the identity in the fusion graph. Indeed, in the conformal case,
multi-interval inclusions are self-dual.

For our results conformal invariance is not necessary, although
conformal nets provide the most interesting situation where they can
be applied. We may deal with an arbitrary net
on $\R$, provided it is strongly additive (a property equivalent to
Haag duality on $\R$ if conformal invariance is assumed) and there
exists a cyclic and separating vector for the von Neumann
algebras of half-lines (vacuum), such that the  corresponding
modular conjugations act geometrically as PCT symmetries
(automatic in the conformal case). We will deal with this more general
context.

Our paper is organized as follows. Then we consider representations. The second section
discusses general properties of multi-interval inclusions and in particular gives
motivations for the strong additivity assumption. The third section enters the core of
our analysis and contains a first inequality between the global index of the sectors and
the index of the 2-interval subfactor. In Section \ref{S4} we study the structure of
sectors associated with the LR net, an analysis mostly based on the braiding symmetry, the
work of Izumi \cite{I4} and the $\a$-induction, which has been introduced in \cite{LR} and
further studied in \cite{X1,BE,BEK}. Section \ref{S5} combines and develops the previous
analysis to obtain our main results for the 2-interval inclusion. These results are
extended to the case of $n$-interval inclusions in Section \ref{S6}. We then we illustrate
our results in models and examples in Section \ref{Examples}. We collect in Appendix
\ref{A1} the results the universal crossed product description of the LR inclusion and of
its multiple iterated occurring in our analysis. We include a further appendix concerning
the disintegration of locally normal or localizable representations into irreducible ones,
that is needed in the paper; these results have however their own interest.

For basic facts concerning conformal nets of von Neumann algebras on 
$\mathbb R$ or $S^1$, the reader is referred to \cite{GL,LR}, see also 
the Appendix \ref{rep-endo}.

\section{General properties}

In this section we shortly examine a few elementary properties for nets 
of von Neumann algebras, partly to motivate our strong additivity 
assumption in the main body of the paper, and partly to examine 
relations with dual nets.  To get our main result, the reader may 
however skip this part, except for Proposition \ref{p9}, and get 
directly to the next section, where we will restrict our study to 
completely rational nets.

In this section, $\A$ will be a local irreducible net of von Neumann
algebras on $S^1$, namely, $\A$ is an inclusion preserving map
$$
\I\ni I\mapsto \A(I)
$$
from the set $\I$ of intervals (open, non-empty sets with
contractible closure)  of $S^1$
to von Neumann algebras on a fixed Hilbert $\H$ space such that
$\A(I_1)$
and $\A(I_2)$ commute if $I_1\cap I_2=\emptyset$ and
$\bigvee_{I\in\I} \A(I)=
B(\H)$, where $\vee$ denotes the von Neumann algebra generated.

If $E\subset S^1$ is any set, we put
$$
\A(E)\equiv \bigvee\{\A(I): \ I\in \I,\ I\subset E\}
$$
and set
$$
\hat \A(E)\equiv \A(E')'
$$
with $E'\equiv S^1\setminus E$.\footnote{The results in this section
are also valid
for nets of von Neumann algebras on $\R$, if $\I$ denotes the set of
non-empty bounded open intervals of $\R$ and $E'=\R\setminus E$ for
$E\subset \R$.}

We shall assume Haag {\it duality} on $S^1$, which automatically
holds if
$\A$ is conformal \cite{BGL}, namely,
$$
\A(I)'=\A(I'), \quad I\in\I,
$$
thus $\hat \A(I)=\A(I)$, $I\in\I$,
but for a disconnected set $E\subset S^1$,
$$
\A(E)\subset \hat \A(E)
$$
is in general a non-trivial inclusion.

We shall say that $E\subset S^1$ is an $n${\it -interval} if both $E$
and $E'$
are unions of $n$ intervals with disjoint closures, namely
$$E=I_1\cup I_2\cup\cdots \cup I_n,\quad I_i\in\I,$$
where $\bar I_i\cap \bar I_j=\emptyset$ if $i\neq j$.
The set of all $n$-intervals will be denoted by $\I_n$.

Recall that $\A$ is $n$-{\it regular},
if $\A(S^1\setminus\{p_1,\dots p_{n}\})=B(\H)$
for any $p_1,\dots p_{n}\in S^1$.

Notice that $\A$ is 2-regular if and only if the $\A(I)$'s are factors,
since we are assuming Haag duality, and that $\A$ is 1-regular if for
each point $p\in S^1$
\begin{equation}\label{point}
\bigcap_n \A(I_n)=\C
\end{equation}
if $I_n\in\I$ and $\bigcap_n I_n = \{p\}$.

\begin{proposition}\label{p1}
The following are equivalent for a fixed $n\in\mathbb N$:
\begin{description}
\item{\rm (i)} The inclusion $\A(E)\subset \hat \A(E)$ is irreducible
for $E\in \I_n$.
\item{\rm (ii)} The net $\A$ is $2n$-regular.
\end{description}
\end{proposition}

\proof With $E=I_1\cup\cdots\cup I_n$ and $p_1,\dots, p_{2n}$ the $2n$
boundary points of $E$, we have $\A(E)'\cap \hat \A(E)=\C$ if and only if
$ \A(E)\vee \hat \A(E)'=B(\H)$, which holds if and  only if
$\A(E)\vee \A(E')=B(\H)$, thus if and only if
$\A(S^1\setminus\{p_1,\dots, p_{2n}\})=B(\H)$,
namely $\A$ is $2n$-regular.
\endproof
If $\A$ is {\it strongly additive}, namely,
$$
\A(I)=\A(I\setminus\{p\})
$$
where $I\in\I$ and $p$ is an interior point of $I$, then $\A$ is
$n$-regular for all $n\in\N$, thus all
$\A(E)\subset\hat \A(E)$ are irreducible inclusions of factors, $E\in\I_n$.

A partial converse holds.

If $\cN\subset\M$ are von Neumann algebras, we shall say that
$\cN\subset \M$ has {\it finite-index} if the Pimsner-Popa inequality
\cite{PP} holds, namely
there exists $\lambda >0$ and a conditional expectation $\E:\M\to\cN$
with
$\E(x)\ge\lambda x$, for all $x\in\M_+$, and denote the index by
$$
[\M:\cN]_\E=\lambda^{-1}
$$
with $\lambda$ the best constant for the inequality to hold and
$$
[\M:\cN]=[\M:\cN]_{\min}=\inf_\E\;[\M:\cN]_\E
$$
denotes the minimal index, (see \cite{K} for an overview).

Recall that $\A$ is {\it split} if there exists an intermediate type I 
factor between $\A(I_1)$ and $\A(I_2)$ whenever $I_1$, $I_2$ are 
intervals and the closure $\bar I_1$ is contained in the interior of 
$I_2$.  This implies (indeed it is equivalent to e.g.  if the 
$\A(I)$'s are factors) that $\A(I_1)\vee \A(I'_2)$ is naturally 
isomorphic to the tensor product of von Neumann algebras 
$\A(I_1)\otimes \A(I'_2)$ (cf.\cite{DL}) .  For a conformal net, the 
split property holds if Tr$(e^{-\beta L_0}) < \infty$ for all 
$\beta>0$, cf.  \cite{DRL}.

Notice that if $\A$ is split and  $\A(I)$ is a factor for $I\in\I$,
then $\A(E)$ is a factor for  $E\in\I_n$ for any $n$.

\begin{proposition}\label{p2}
Let $\A$ be split and 1-regular. If there exists a constant $C>0$ such
that
$$
[\hat \A(E):\A(E)]<C \quad \forall\quad E\in \I_2,
$$
then
$$
[\A(I):\A(I\setminus\{p\})]<C \quad\forall I\in\I,\ p\in I.
$$
\end{proposition}

\proof With $I\in\I_2$ and $p\in I$ an interior point, let $I_1,
I_2\in\I$ be
the connected components of $I\setminus\{p\}$, let $I^{(n)}_2\subset
I_2$ be an increasing  sequence of intervals with one boundary point
in
common with $I$ such that $p\notin \ov{I_2^{(n)}}$ and
$\bigcup_nI_2^{(n)}=I_2$.
Then $E_n\equiv I_1\cup I_2^{(n)}\in \I_2$ and we have
\begin{eqnarray*}
&&\A(E_n)\nearrow \A(I\setminus\{p\}),\\
&&\hat \A(E_n)\nearrow \A(I),
\end{eqnarray*}
where $\cN_n\nearrow\cN$ means $\cN_1\subset\cN_2\subset\cdots$ and
$\cN=\bigvee \cN_n$, while $\cN_n\searrow\cN$ will mean
$\cN_1\supset\cN_2
\supset\cdots$ and
$\cN=\bigcap \cN_n$. The first relation is clear by definition.  The
second
relation follows because
$$
\hat \A(E_n)'=\A(E'_n)=\A(I')\vee \A(L_n),
$$
where $E'_n\in\I_2$, $E_n=I'\cup L_n$, and $\bigcap L_n=\{p\}$,
therefore
$\A(L_n)\searrow {\cal C}$. By the split property $\A(I')\vee \A(L_n)\isom
\A(I')\otimes \A(L_n)$, hence by eq. (\ref{point})
$$
\A(E'_n)\searrow \A(I'),
$$
thus
$$
\hat \A(E_n)\nearrow \A(I).
$$
The rest of the proof is the consequence of the following general
proposition.
\endproof

\begin{proposition} \label{p3}
\begin{description}
\item{\rm a)} Let
$$\begin{array}{ccccccc}
\cN_1&\subset&\cN_2&\subset&\cdots&\subset& \cN\\
\cap&& \cap &&&& \cap\\
\M_1&\subset&\M_2&\subset&\cdots&\subset& \M
\end{array}$$
be von Neumann algebras, $\cN=\bigvee \cN_i$, $\M=\bigvee \M_i$,
\item{\rm b)} or let
$$\begin{array}{ccccccc}
\cN_1&\supset&\cN_2&\supset&\cdots&\supset& \cN\\
\cap&& \cap &&&& \cap\\
\M_1&\supset&\M_2&\supset&\cdots&\supset& \M
\end{array}$$
be von Neumann algebras, $\cN=\bigcap \cN_i$, $\M=\bigcap \M_i$.
\end{description}

\noindent Then
$$[\M:\cN]\le \liminf_{i\to\infty}\;[\M_i:\cN_i].
$$
\end{proposition}
\proof It is sufficient to prove the result in the situation b) as
the case
a) will follow after taking commutants.
We may assume $\liminf_{i\to\infty}[\M_i:\cN_i]<\infty$.

Let $\E_i:\M_i\to\cN_i$ be an expectation and
$\lambda>\liminf_{i\to\infty}[\M_i,:\cN_i]_{\E_i}$. Then there exists
$i_0$
such that
for all $x\in\M_i^{+}$, $i\ge i_0$,
$$
\E_i(x)\ge\lambda^{-1}x.
$$
Let $\E^{(0)}_i=\E_i\vert_\M$, considered as a map from $\M$ to
$\cN_i$,
and let $\E$ be a weak limit point of $\E_i^{(0)}$. Then
$$
\E(x)\ge\lambda^{-1}x,\qquad x\in\M_+,
$$
and $\E(\M)\subset \bigcap_i \cN_i=\cN$, moreover
$\E\vert_\cN=\hbox{\rm
id}$, because $\E_i\vert_\cN=\hbox{\rm id}$. Thus $\E$ is an
expectation of $\M$ onto $\cN$ and
$$
[\M:\cN]\le [\M:\cN]_{\E}\le\lambda.
$$
As $\E_i$ is arbitrary, we thus have $[\M:\cN]\le
\liminf_{i\to\infty}[\M_i,:\cN_i]$.
\endproof

Recall now that the {\it dual net} $\A^d$ of $\A$ is the net on the
intervals of $\R$ defined by $\A^d(I)\equiv \A(\R\setminus I)'$, where
we have chosen a point $\infty\in S^1$ and identified $S^1$ with
$\R\cup\{\infty\}$.

Note that if $\A$ is conformal, then Haag duality
automatically holds \cite{GLW} and the dual net $\A^d$ is also
a conformal net which is moreover strongly additive;
furthermore $\A=\A^d$, if and only if $\A$
is strongly additive, if and only if Haag duality holds on $\R$.

\begin{corollary}\label{c8}
In the hypothesis of Proposition \ref{p2}, let $\A^d$ be
the dual net on $\R$, then
$$
\A(I)\subset \A^d(I)
$$has finite index for all bounded intervals $I$ of $\R$.
\end{corollary}

\proof Denoting $I_1 = I'$, the complement of $I$ in $S^1$,
the commutant of the inclusion $\A(I)\subset \A^d(I)$ is
$\A(I_1\setminus\{\infty\})\subset \A(I_1)$, and this has finite
index.
\endproof

We have no example where $\A(I)\subset \A^d(I)$ is non-trivial with
finite
index and $\A$ is conformal; therefore the equality $\A(I)=\A^d(I)$,
i.e.
strong additivity, might
follow from the assumptions in Corollary \ref{p2} in the conformal
case.

\begin{proposition}\label{p9}
Let $\A$ be split and strongly additive, then
\begin{description}
\item{\rm (a)} The index $[\hat \A(E):\A(E)]$ is independent of
$E\in\I_2$.
\item{\rm (b)} The inclusion $\A(E)\subset\hat \A(E)$ is irreducible
for $E\in\I_2$.
\end{description}
\end{proposition}

\proof Statement (b) is immediate by Proposition \ref{p1}.

Concerning (a), let
$E=I_1\cup I_2$ and $\tilde E=I_1\cup \tilde I_2$ where $\tilde
I_2\supset
I_2$ are intervals and $I_0\equiv\tilde I_2\setminus I_2$. Assuming
$\lambda^{-1}\equiv[\hat \A(\tilde E):\A(\tilde E)]<\infty$, let
$\E_{\tilde E}$
be the corresponding expectation with $\lambda$-bound. Of course
$\E_{\tilde E}$ is the identity on $\A(I_0)$, hence
$$
\E_{\tilde E}(\hat \A(E))\subset \A(I_0)'\cap \A(\tilde E)=\A(E)
$$
where last equality follows at once by the split property and strong
additivity as $\A(I_0)'\cap \A(\tilde I_2)=\A(I_2)$.

Therefore $\E_{\tilde E}\mid_{\hat \A(E)}=\E_E$ showing
$$
[\hat \A(E):\A(E)]\le [\hat \A(\tilde E):\A(\tilde E)],
$$
where we omit the symbol ``min'' as the expectation is unique. Thus
the index
decreases by decreasing the $2$-interval. Taking commutants, it also
increases, hence it is constant.
\endproof

\begin{corollary}\label{c10}
Let $\A$ satisfy the assumption of Proposition \ref{p2} and let
$\A^d$ be
the dual net on $\R$ of $\A$.
Then
$$
[\widehat{\A^d}(E):\A^d(E)]<\infty \quad\forall E\in \I_2.
$$
\end{corollary}

\proof We fix the point $\infty$ and may assume $E=I_1\cup I_2$ with
$\infty\in I_2$. Set $E'=I_3\cup I_4$ with $I_3\ni \infty$. Then
$\A^d(I_3)=\A(I_3)$, $\A^d(I_2)=\A(I_2)$ and we have
\begin{eqnarray*}
\A(E)&\subset& \A^d(I_1)\vee \A(I_2) \\
&=&\A^d(E)\subset\widehat{\A^d}(E)\\
&=&(\A(I_3)\vee \A^d(I_4))'
\subset(\A(I_3)\vee \A(I_4))'=\hat \A(E).
\end{eqnarray*}
\endproof
Anticipating results in the following, we have:
\begin{corollary}\label{n-reg}
Let $\A$ be a local irreducible conformal split net on $S^{1}$. If
$[\hat \A(E):\A(E)]=\In_{\textnormal{global}}
<\infty$, $E\in\I_2$, then $\A$ is $n$-regular
for all $n\in\mathbb N$.
\end{corollary}
\proof
If $\rho$ is an irreducible endomorphism of $\A$ localized
in an interval $I$, then $\rho |_{\A(I)}$ is irreducible \cite{GL}.
Therefore, by Th. \ref{t11} (and comments there after) and Prop. \ref{p22}, the
assumptions
imply that if $E\in\I_2$ then $\A(E)\subset\hat \A(E)$ is the LR
inclusion associated with the system of all irreducible sectors,
which is irreducible. Then $\A(E)\subset\hat \A(E)$
is irreducible for all $E\in\I_n$ as we shall see in Sect. \ref{S6}.
By Prop. \ref{p1} this implies the regularity for all $n$.
\endproof

In view of the above results,
it is natural to deal with strongly additive nets,
when considering multi-interval inclusions of local algebras and thus
to deal with nets of factors on $\R$, as we shall do in the
following.

\section{Completely rational nets}
In this section we will introduce the notion of completely rational 
net, that will be the main object of our study in this paper, and 
get a first analysis.

In the following, we shall denote by $\I$ the set of bounded open 
non-empty intervals of $\R$, set $I'=\R\setminus I$ and define 
$\A(E)=\bigvee \{\A(I),I\subset E, I\in\I\}$ for $E\subset\R$. 
We again denote by $\I_n$ the set of unions of  $n$ elements of 
$\I$ with pairwise disjoint closures.  
\footnote{There will be no conflict with the notations in the 
previous section as the point $\infty$ does not contribute to the 
local algebras and we may extend $\A$ to $S^1$ setting $\A(I)\equiv 
\A(I\setminus\{\infty\})$, see Appendix \ref{rep-endo}.}

\begin{definition}
A local irreducible net $\A$ of von Neumann algebras on the intervals of 
$\R$ is called {\it completely rational} if the following holds:
\begin{description}
\item{\rm (a)} Haag duality on $\R: \ \A(I')=\A(I)'$, $I\in \I$,
\item{\rm (b)} $\A$ is strongly additive,
\item{\rm (c)} $\A$ satisfies the split property,
\item{\rm (d)}  $[\hat \A(E):\A(E)]<\infty$, if $E\in\I_2$,
\end{description}
\end{definition}

Note that, if $\A$ is the restriction to $\mathbb R$ of a 
local conformal net on $S^1$ (namely a local net which is M\"obius 
covariant with positive energy and cyclic vacuum vector) then (a) is 
equivalent to (b), cf. \cite{GLW}.

We shall denote by $\mu_\A=[\hat \A(E):\A(E)]$ the index of the
irreducible
inclusion of factors $\A(E)\subset \hat \A(E)$ in case $\mu_\A$ is
independent
of $E\in\I_2$, in particular if $\A$ is split, by Proposition
\ref{p9}.

By a {\it sector} $[\rho]$ of $\A$ we shall mean the equivalence class 
of a localized endomorphism $\rho$ of $\A$, that will always be 
assumed to be transportable i.e.  localizable in each bounded 
interval $I$ (see also Appendix \ref{rep-endo}).  Unless otherwise 
specified, a localized endomorphism $\rho$ has finite dimension.  If 
$\rho$ is localized in the interval $I$, its restriction 
$\rho|_{\A(I)}$ is an endomorphism of $\A(I)$, thus it gives rise to a 
{\it sector} of the factor $\A(I)$ (i.e.  a normal unital endomorphism 
of $\A(I)$ modulo inner automorphisms of $\A(I)$ \cite{L2}) and it 
will be clear from the context which sense will be attributed to the 
term sector.

The reader unfamiliar with the sector strucure is referres to 
\cite{L2,LR,GL} and to the Appendix \ref{rep-endo}.

Let $E=I_1\cup I_2\in \I_2$ and $\rho$ and $\sigma$ irreducible
endomorphisms of $\A$ localized respectively in $I_1$ and in $I_2$.
Then $\rho\sigma$ restricts to an endomorphism of
$\A(E)$, since both $\rho$ and $\sigma$  restrict.

Denote by $\gamma_E$ the canonical endomorphism of $\hat \A(E)$ into
$\A(E)$ and $\lambda_E\equiv\gamma_E\vert_{\A(E)}$.

\begin{theorem}\label{t11} Let $\A$ be  completely rational. With the above
notations,
$\rho\sigma\vert_{\A(E)}$ is contained in $\lambda_E$ if and only if
$\sigma$ is
conjugate to $\rho$. In this case
$\rho\sigma\vert_{\A(E)}\prec\lambda_E$
 with multiplicity one.\end{theorem}
\proof
 By \cite{LR} $\rho\sigma |_{\A(E)}\prec
\lambda_E$ if and only if there exists
an isometry $v\in \hat \A(E)$ such that
\begin{equation}
        vx=\rho\sigma(x)v \quad\forall x\in \A(E).
        \label{uni}
\end{equation}
If equation (\ref{uni}) holds, then it holds for $x\in\A(I)$ for
all $I\in\I$ by strong additivity, hence $\sigma=\bar\rho$.

Conversely, if $\sigma=\bar\rho$, then there exists an isometry $v\in
\A(I)$ such that $vx=\rho\sigma(x)v$ for all $x\in \A(I)$, where
$I$ is the interval  $I\supset E$ given by $I=I_1\cup I_2\cup \bar
I_3$
with $I_3$
the bounded connected component of $E'$.

Since $\rho$ and $\sigma$ act trivially on $\A(I_3)$, we have
$$
v\in \A(I_3)'\cap \A(I),
$$
but
$$
\A(I_3)'\cap \A(I)=(\A(I_3)\vee \A(I'))'=\A(E')'=\hat \A(E),
$$
therefore equation (\ref{uni}) holds true.
As the $\rho$ and $\sigma$ are irreducible, the isometry $v$ in
(\ref{uni})
unique up to a phase and this is equivalent to
$\rho\bar\rho\vert_{\A(E)}\prec\lambda_E$ with multiplicity one.
\endproof

We remark that in the above theorem strong additivity is not necessary
 for $\rho\bar\rho\prec \lambda_E$,
as can be replaced by the factoriality of $\A(E)$, equivalently
of $\hat \A(E)$; this holds e.g. in the conformal case.

Moreover also the split property is unnecessary, it has not been used.

We shall say that the net $\A$ on $\R$ has a {\it modular PCT}
symmetry, if there exists a cyclic separating (vacuum) vector $\Omega$
for each
$\A(I)$, if $I$ is a half-line (Reeh-Schlieder property), and the
modular conjugation $J$ of $\A(a,\infty)$ with
respect to $\Omega$ has the geometric property
\begin{equation}
J\A(I + a)J=\A(-I+a),\quad I\in\I,\ \ a\in \R.
	\label{pct}
\end{equation}
This is automatic if $\A$ is conformal, see \cite{BGL,FG}. It easy to see that
the modular PCT property implies translation covariance, where the
translation unitaries are products of modular
conjugations, but positivity of the energy does not necessarily holds.

Note that eq. (\ref{pct}) implies Haag duality for half-lines
\[
\A(-\infty,a)'=\A(a,\infty),\quad a\in\mathbb R\ .
\]
Setting $j\equiv{\rm \Ad}J$, the conjugate sector exists and it is
given
by the formula \cite{GL1}
$$
\bar\rho=j\cdot\rho\cdot j.
$$
\begin{corollary}\label{c12}
If $\A$ is completely rational with modular PCT, then $\A$ is
rational,
namely there are only finitely many irreducible sectors
$[\rho_0], [\rho_1], \dots, [\rho_n]$ with finite dimension and we
have
\begin{equation}
\sum^n_{i=0}d(\rho_i)^2\le\mu_\A.
\label{ineq}\end{equation}
\end{corollary}

\proof It is sufficient to show this last inequality. By the split
property, the endomorphisms $\rho_i\bar\rho_i\vert_{\A(E)}$ can be
identified with the endomorphisms $\rho_i\otimes\bar\rho_i$ on
$\A(I_1)\otimes \A(I_2)$, hence they are mutually inequivalent.

By Theorem \ref{t11},
\begin{equation}
 \bigoplus^n_{i=1}\rho_i\bar\rho_i\vert_{\A(E)} \prec \lambda_E,
\label{rho_2}\end{equation}
hence
$$
\mu_\A=[\hat \A(E):\A(E)]=d(\lambda_E)\ge\sum d(\rho_i)^2.
$$
\endproof

We now give a partial converse to Theorem \ref{t11}.

\begin{lemma} \label{l14} Let $\A$ be completely rational  and let
$\E_E$ be the conditional
expectation $\hat \A(E)\to \A(E)$.
\begin{description}
\item{\rm (a)} If $E\subset\tilde E$ and $E,\tilde E\in \I_2$, then
$$
\E_{\tilde
E}\vert_{\hat \A(E)}=\E_E.
$$
\item{\rm (b)} There exists a canonical endomorphism $\gamma_{\tilde
E}$
of $\hat \A(\tilde E)$ to $\A(\tilde E)$
such that $\gamma\vert_{\hat \A(E)}$ is a
canonical endomorphism of $\hat \A(E)$ into $\A(E)$ and satisfies
$$
\gamma\vert_{\hat \A(E)'\cap \A(\tilde E)}=\hbox{\rm id}.
$$\end{description}
\end{lemma}

\proof
(a) has been shown in the proof of Proposition \ref{p9}.

(b) is an immediate variation of \cite[Proposition 2.3]{GL1} and
\cite[Theorem 3.2]{LR}.
\endproof
\begin{theorem}\label{t13b} Let $\A$ be completely rational.
Given $E\in
\I_2$, $\lambda_E$ extends to a localized (transportable)
endomorphism $\lambda$ of
$\A$ such that $\lambda\vert_{\A(I)}=\hbox{\rm id}$, if $I\subset E'$,
$I\in\I$. Moreover, $d(\lambda)=d(\lambda_E)=\mu_\A$.

In particular, if $\A$ is conformal, then $\lambda$ is M\"obius
covariant with positive energy.\end{theorem}

\proof
Let $E=(a,b)\cup(c,d)$ where $a<b<c<d$ and $\tilde E=(a',b)\cup(c,d')$
where $a'<a$ and $d'>d$. By Lemma \ref{l14} we have a $\gamma_{\tilde
E}$ with
$\lambda_{\tilde E}\vert_{\A(I)}=\hbox{\rm id}$, if $I\subset \I$,
$I\in\tilde E\setminus E$.

Analogously there is a canonical endomorphism
$\gamma:\hat \A(\tilde E)\to \A(\tilde E)$ acting trivially
on $\A(E)$. We may write
$$\gamma_{\tilde E}=\Ad\ u\cdot\gamma
$$
with $u\in \A(\tilde E)$, hence
$$
\lambda_{\tilde E}=\Ad u\cdot \lambda,\quad
\lambda=\gamma\vert_{\A(\tilde E)}.
$$
Since $\gamma\mid_{\A(a,b)}=\hbox{\rm id}$, $\gamma\mid_
{\A(c,d)}=\hbox{\rm id}$, we have
$$
\lambda_{\tilde E}=\Ad u\quad \hbox{\rm on\ }
\A(a,b), \A(c,d).
$$
Therefore, the formula
$$\tilde\lambda=\Ad u $$
defines an endomorphism of $\A(a,d)$ acting  trivially
an $\A(b,c)$, with
$$\tilde \lambda\vert_{\A((a,b)\cup (c,d))}=\lambda_E.
$$
We may also have chosen $\gamma$ ``localized''  in $(a',a'')\cup
(d'',d')$
with $a'<a''<a$ and $d<d''<d'$
so that we may assume $\tilde\lambda$ to act trivially on
$\A((a'',b)\vee (c,d'')).$

Letting $a',a''\to -\infty$ and $d'',d'\to +\infty$,
we construct, by an inductive limit of the
$\tilde \lambda$'s, an endomorphism $\lambda$ of
the quasi-local $C^*$-algebra $\overline{\bigcup_{s>0} \A(-s,s)}$.

Clearly, $\lambda$ is localized in $(a,d)$,
acts trivially on $\A(b,c)$ and is transportable. Moreover, $\lambda$
has finite index as the operators $R,\bar R\in(i,\lambda^2)$
in the standard solution for the conjugate equation
\cite{L2,LRo}
$$
\bar R^*\bar\lambda (R)=1,\quad R^*\lambda (\bar R)=1,
$$
on $\hat \A(E)$ give the same relation on $\A(I)$ for any $
I\supset E$, $I\in\I$.

If $\A$ is conformal, then $\rho$ is covariant with
respect to translations and dilations by \cite{GL}.
As we may vary the point $\infty$, $\lambda$ is covariant with respect
to dilations and translations with respect to different point at
$\infty$,
hence $\lambda$ is M\"obius covariant.
\endproof
\begin{lemma}\label{L12} Let $\A$ be completely rational.  Then there 
are at most $\lfloor\mu_\A\rfloor$ mutually different irreducible 
sectors of $\A$ (with finite or infinite dimension).  
\end{lemma} 
\proof 
Consider the family $\{[\rho_{\lambda}]\}$ of all 
irreducible sectors and let $N$ be the cardinality of this family.  
With $E=I_1\cup I_2 \in\I_2$, we may assume that each $\rho_{\lambda}$ 
is localized in $I_1$ and choose endomorphisms $\sigma_{\lambda}$ 
equivalent to $\rho_{\lambda}$ and localized in $I_2$.  Let then 
$u_{\lambda}\in(\rho_{\lambda},\sigma_{\lambda})\subset \hat \A(E)$ be 
a unitary intertwiner and $\E$ the conditional expectation from $\hat 
\A(E)$ to $\A(E)$.  Since \[ u_{\lambda}\rho_{\lambda}(x)= 
\sigma_{\lambda}(x)u_{\lambda} =xu_{\lambda}\ , \quad \forall x\in 
\A(I_1), \] we have \[ u^*_{\lambda'}u_{\lambda}\,\rho_{\lambda}(x)= 
\rho_{\lambda'}(x)\,u^*_{\lambda'}u_{\lambda}\ , \quad\forall x\in 
\A(I_1), \] hence $T=\E(u^*_{\lambda'}u_{\lambda})\in\2A(E)$ 
intertwines $\rho_{\lambda}\vert\A(I_1)$ and 
$\rho_{\lambda'}\vert\A(I_1)$.  The split property allowing us to 
identify $\A(E)$ and $\A(I_1)\otimes\A(I_2)$, every state $\phi$ in 
$\A(I_2)_*$ gives rise to a conditional expectation 
$\E_\phi:\A(E)\rightarrow\A(I_1)$.  Then 
$\E_\phi(T)\in(\rho_{\lambda},\rho_{\lambda'})$, and the inequivalence 
of $\rho_\lambda\vert\A(I_1), \rho_{\lambda'}\vert\A(I_1)$, see above, 
entails $\E_\phi(T)=0$.  Since this holds for every $\phi\in\A(I_2)_*$ 
we conclude \[ T=\E(u^*_{\lambda'}u_{\lambda})=0, \quad \lambda'\neq 
\lambda.  \] Let $\M$ be the Jones extension of $\A(E)\subset \hat 
\A(E)$, $e\in \M$ the Jones projection implementing $\E$ and let 
$\2E_1:\2M\rightarrow\hat{\2A}(E)$ be the dual conditional 
expectation.  Then $eu^*_{\lambda'}u_{\lambda}e =0$ if $\lambda'\neq 
\lambda$ and therefore the $e_{\lambda}\equiv u_{\lambda}e 
u^*_{\lambda}$ are mutually orthogonal projections in $\2M$ with 
$\2E_1(e_{\lambda})=\mu_\A^{-1}$.  Since their (strong) sum 
$p=\sum_\lambda e_\lambda$ is again an orthogonal projection we have 
$p\le 1$ and thus $\2E_1(p)\le\2E(1)=1$.  This implies the bound 
$N\mu_\A^{-1}\le 1$ and thus our claim.  
\endproof
We shall say that a sector $[\rho]$ is of {\it type I} if 
$\vee_{I\in\I}\rho(\A(I))$ is a type I von Neumann algebra, namely 
$\rho$ is a type I representation of the quasi local C$^*$-algebra 
$\overline{\cup_{s>o}\A(-s,s)}$.  
\begin{corollary}\label{type I} 
If $\A$ is completely rational on a separable Hilbert space, then all 
factor representations of $\A$ on separable Hilbert spaces are of type 
I. \end{corollary} 
\proof 
Assuming the contrary, by Corollary \ref{inf} we get an infinite 
family $[\rho_{\lambda}]$ of different irreducible sectors.  This is 
in contradiction with the preceding proposition. 
\endproof 
We end this section with the following variation of a known 
fact \cite{DL}.
\begin{proposition}\label{separable}
Let $\A$ be a completely rational net with 
modular PCT on a Hilbert space $\H$. Then $\H$ is separable.
\end{proposition}
\proof
We chose a pair $I\subset\tilde I$ of intervals and a type I factor 
$\N$ between $\A(I)$ and $\A(\tilde I)$.  The vacuum vector $\Omega$ 
is separating for $\A(\tilde I)$, hence for $\N$.  Thus $\N$ admits a 
faithful normal state, hence it is countably decomposable.  Being of 
type I, $\N$ is countably generated.  So 
$\overline{\A(I)\Omega}\subset\overline{\N\Omega}$ is a separable 
subspace of $\H$.  But $\cup_{n=1}^{\infty}\A(-n,n)\Omega$ is dense in 
$\H$, thus $\H$ is separable.
\endproof
\section{The structure of sectors for the (time $= 0$) LR net}
\label{S4} 
This section contains a study of the sector strucure for the net 
obtained by the LR construction, by means of the braiding symmetry. It 
will be continued in the next section by a different approach.

Let $\N$ be an infinite factor and $\{[\rho_i]\}$ a 
rational system of sectors of $\N$, namely the $[\rho_i]$'s form a 
finite family of mutually different irreducible finite-dimensional 
sectors of $\N$ which is closed under conjugation and taking the 
irreducible components of compositions. The identity sector is usually 
labeled as $\rho_0$. We call
$$
\M\supset \N\otimes \N^\o
$$
the {\it LR inclusion}, the canonical inclusion constructed in
\cite{LR} where $\M$ is a factor, $\N\otimes \N^\o\subset \M$ is
irreducible with finite index and
$$
\lambda=\bigoplus_i \rho_i\otimes \rho_i^\o
$$
for $\lambda\in {\rm End}(\N\otimes \N^\o)$ as
the restriction of $\gamma:\M\to \N\otimes \N^\o$. We shall give an
alternative characterization of this inclusion in Proposition
\ref{characterization}.

The same construction works in slightly more generality,
by replacing $\N^\o$ with a factor $\N_1$ and
$\{\rho_i^\o\}_i$ by $\{\rho^j_i\}_i\subset\textnormal{End}(\N_1)$
where $\rho\to\rho^j$ is an
anti-linear invertible tensor functor of the tensor category
generated by $\{\rho_i\}_i$ to the tensor category generated by
$\{\rho_i^j\}_i$. Extensions of our results to this case are obvious,
but sometimes useful, and will be considered possibly implicitly.

The following is due to Izumi \cite{I4}. Since it is easy
to give a proof in our context, we include a proof here.

\begin{lemma}
\label{LL11}
For every $\rho_i$, the $(\N\otimes\N^\o)$-$\M$ sector
$[\rho_i\otimes\id][\ga]
=[\id\otimes\bar\rho_i^\o][\ga]$ is irreducible
and each irreducible
$(\N\otimes\N^\o)$-$\M$ sector arising from $\N\otimes\N^\o\subset
\M$
is of this form,
where $\ga$ is regarded as an $(\N\otimes\N^\o)$-$\M$ sector.
If $[\rho_i]\neq[\rho_j]$ as $\A$-$\A$ sectors, then
$[\rho_i\otimes\id][\ga]\neq [\rho_j\otimes\id][\ga]$ as
$(\N\otimes\N^\o)$-$\M$ sectors.
We have $[\rho_i\otimes\rho_j^\o][\ga]=
\sum_k N_{i\bar j}^k [\rho_k\otimes\id][\ga]$ as
$(\N\otimes\N^\o)$-$\M$ sectors,
where $N_{i\bar j}^k$ is the structure constant for $\{\rho_i\}_i$.
\end{lemma}

\begin{proof}
Set $[\si]=[\rho_i\otimes\id][\ga]$ and compute $[\si][\bar \si]$.
Since $[\bar\ga]=[\iota]$, where $\iota$ is the inclusion map
of $\N\otimes\N^\o$ into $\M$ regarded as a $\M$-$(\N\otimes\N^\o)$
sector, and $[\ga][\iota]=[\lambda]=\sum_k [\rho_k\otimes\rho^\o_k]$, we
have
$[\si][\bar\si]=\sum_k [\rho_i \rho_k \bar\rho_i\otimes\rho_k^\o]$,
and this contains the identity only once.
So $[\rho_i\otimes\id][\ga]$  is an irreducible
$(\N\otimes\N^\o)$-$\M$ sector.
We can similarly prove that if $[\rho_i]\neq[\rho_j]$, then
$[\rho_i\otimes\id][\ga]\neq [\rho_j\otimes\id][\ga]$.

We next set $[\si']=[\id\otimes\bar\rho_i^\o][\ga]$  as
an $(\N\otimes\N^\o)$-$\M$ sector, which is also irreducible.  We
compute
$$[\si][\bar\si']=
[\rho_i\otimes\id][\lambda] [\id\otimes\rho_i^\o]=
\sum_k [\rho_i \rho_k \otimes\rho_k^\o\rho_i^\o],$$
which contains the identity only once.
So we have $[\rho_i\otimes\id][\ga]
=[\id\otimes\bar\rho_i^\o][\ga]$.

The rest is now easy.
\end{proof}

Let us now assume we have
a  strongly additive, Haag dual,  irreducible
net of factors $\A(I)$ on $\R$ with  a rational {\it system of
irreducible sectors} $\{[\rho_i]\}_i$ (with $\rho_0=\id$),
namely $\{[\rho_i]\}_i$ is a family of finitely many
different irreducible sectors of $\A$ with finite dimension
stable under conjugation and irreducible components of compositions.

One may construct \cite{R2,LR} a net of subfactors $\AA\subset \B$ so 
that the corresponding canonical endomorphism restricted on $\AA$ is 
given by $\bigoplus_i \rho_i\otimes\rho_i^\o$.  We call this $\B$ the 
{\it LR net}.  For $\A^\o$, we use 
$\e^\o(\rho_k^\o,\rho_l^\o)=j(\e(\rho_k,\rho_l))^*$, where $j$ is the 
anti-isomorphism from $\A$ to $\A^\o$.  In order to distinguish two 
braidings, we write $\e^+$ and $\e^-$.

In other words, the LR net here is obtained as the time zero fields 
from the canonical two-dimensional net constructed in \cite{LR}: it is 
a local net, but if $\A$ is translation covariant with positive 
energy, $\B$ is translation covariant without the spectrum condition 
(the translation on $\B$ are space translations).

Then the net of inclusion $\AA(I)\subset \B(I)$ is a net of subfactors 
in the sense of \cite[Section 3]{LR}, that is, we have a vacuum vector 
with Reeh-Schlieder property and consistent conditional expectations.  
We denote by $\ga$ the canonical endomorphism of $\B$ into $\AA$ and 
its restriction to $\AA$ by $\lambda$.  We may suppose that also 
$\lambda$ is localized in $I$.  We shorten our notation by setting 
$\N\equiv\A(I)$ and $\M=\B(I)$.  We thus have $\lambda(x)=\sum_i 
V_i(\rho_i\otimes \rho^\o_i)(x)V_i^*$, where $V_i$'s are isometries in 
$\N\otimes\N^\o$ with $\sum_i V_i V_i^*=1$.

We follow \cite{I3} for the terminology of $(\N\otimes\N^\o)$-$\M$ 
sectors, and so on, and study the sector structure of the subfactor 
$\N\otimes\N^\o\subset \M$ in this section.  In other words we study 
the sector structure of a single subfactor, not the structure of 
superselection sectors of the net, though we will be interested in 
this structure for the net in the next section.  So the terminology 
{\it sector} is used for a subfactor, not for a net, in this section.  
However the inclusion $\N\otimes\N^\o\subset\M$ has extra structure 
inherited by the inclusion of nets $\AA\subset\B$, that is there are 
the left and right unitary braid symmetries and the extension and 
restriction maps.  We first note that 
$\{[\rho_i\otimes\rho_j^\o]\}_{ij}$ gives a system of irreducible 
$\AA$-$\AA$ sectors.

This gives the description of the principal graph of 
$\N\otimes\N^\o\subset \M$ as a corollary as follows, which was first 
found by Ocneanu in \cite{O2} for his asymptotic inclusion.  Label 
even vertices with $(i,j)$ for $[\rho_i\otimes\rho_{\bar j}^\o]$ and 
odd vertices with $k$ for $[\rho_k\otimes\id][\ga]$ and draw an edge 
with multiplicity $N_{ij}^k$ between the even vertex $(i,j)$ and the 
odd vertex $k$.  The connected component of this graph containing the 
vertex $(0,0)$ is the principal graph of the subfactor 
$\N\otimes\N^\o\subset \M$.

Now we consider the $\a$-induction introduced in \cite{LR} and further 
studied in \cite{X1, BE}, namely if $\sigma$ is a localized 
endomorphism of $\AA$, we set 
\begin{equation} 
\a^\pm_{\sigma}= 
\ga^{-1}\cdot\Ad(\e^\pm(\sigma,\lambda))\cdot\sigma\cdot\ga.  
\label{induction} 
\end{equation} 
(The notation in \cite{LR} is $\sigma^{\ext}$).

Recall that if $\sigma$ is an endomorphism of $\AA$
localized in the interval $I$, then
$\a^\pm_{\sigma}$ is an endomorphism of $\B$ localized
in a positive/negative half-line containing
$I$, yet, as shown in \cite[I]{BE}, $\a^\pm_{\sigma}$ restricts to an
endomorphism of $\M=\B(I)$. We will denote this restriction
by the same symbol $\a^\pm_{\sigma}$.
\begin{lemma}
\label{LL12}
The $\M$-$\M$ sectors $[\a^+_{\rho_i\otimes\id}]$ are irreducible
and  mutually different.
\end{lemma}

\begin{proof}
We compute
$\lan \a^+_{\rho_i\otimes\id}, \a^+_{\rho_j\otimes\id}\ran$, the
dimension of the intertwiner space between
$\a^+_{\rho_i\otimes\id}$
and $\a^+_{\rho_j\otimes\id}$, by using
\cite[I, Theorem 3.9]{BE}.  This number is then equal to
$$\lan \bigoplus_k \rho_k\rho_i\otimes\rho_k^\o,\rho_j\otimes\id\ran
=\de_{ij}.$$
This gives the conclusion.
\end{proof}

\begin{lemma}
\label{LL13}
As $\M$-$\M$ sectors, we have
$[\a^+_{\rho_i\otimes\id}]
=[\a^+_{\id\otimes\rho_i^\o}]$.
\end{lemma}

\begin{proof}
By a similar argument to the proof of the above lemma, we know
that $[\a^+_{\id\otimes\rho_i^\o}]$ is also irreducible.
\cite[I, Theorem 3.9]{BE} gives
$$\lan [\a^+_{\rho_i\otimes\id}],
[\a^+_{\id\otimes\rho_i^\o}]\ran=
\lan \bigoplus_k \rho_k\rho_i\otimes\rho_k^\o,\id\otimes\rho_i^\o\ran
=1,$$
which gives the conclusion.
\end{proof}

We then have the following corollary.

\begin{corollary}
\label{CC14}
The set of irreducible $\M$-$\M$ sectors appearing in the
decomposition of $\a^+_{\rho_i\otimes\rho_j^\o}$ for all $i,j$
is $\{[\a^+_{\rho_i\otimes\id}]\}_i$.
\end{corollary}

The next theorem is useful for studying the
subfactors arising from disconnected intervals for a conformal
net. For the rest of this section we shall assume the
braiding to be non-degenerate.

\begin{theorem}
\label{intersect}
Assume the braiding to be non-degenerate and
suppose an irreducible $\M$-$\M$ sector $[\be]$ appears in
decompositions
of both $\a^+_{\rho_i\otimes\rho_j^\o}$
and $\a^-_{\rho_k\otimes\rho_l^\o}$  for some $i,j, k, l$.
Then $[\be]$ is the identity of $\M$.
\end{theorem}

\begin{proof}
$\a^+$ and $\a^-$ map sectors localized in
bounded intervals to soliton sectors localized in right unbounded
and left unbounded half-lines, respectively. Hence $[\be]$ is
localized
in a bounded interval.
By the above corollary, we may assume that
$[\be]=[\a^+_{\rho_i\otimes\id}]$
for some $i$, hence $\rho_i\otimes\id$ must have trivial monodromy
with $\lambda$, i.e.,
$\e(\rho_i\otimes\id,\lambda)\e(\lambda,\rho_i\otimes\id)=1$,
which in turn gives  $\e(\rho_i,\rho_k)\e(\rho_k,\rho_i)=1$
for all $k$.  The non-degeneracy assumption gives
$[\rho_i]=[\id]$ as desired.
\end{proof}

We now define an endomorphism of $\M$ by
$\be_{ij}=\a^+_{\rho_i\otimes\id}\a^-_{\id\otimes\rho_j^\o}$.
More explicitly, we have
$\be_{ij}=\ga^{-1}\cdot\Ad(U_{ij}^{+-})\cdot(\rho_i\otimes
\rho_j^\o)\cdot\ga$, where
$$U_{ij}^{+-}=\sum_k V_k(\e^+(\rho_i,\rho_k)\otimes
\e^{-,\o}(\rho_j^\o,\rho_k^\o))(\rho_i\otimes\rho_j^\o)(V_k^*).$$
Note that if we define similarly
$$U_{ij}^{++}=\sum_k V_k(\e^+(\rho_i,\rho_k)\otimes
\e^{+,\o}(\rho_j^\o,\rho_k^\o))(\rho_i\otimes\rho_j^\o)(V_k^*),$$
we then have $\a^+_{\rho_i\otimes\rho_j^\o}=
\ga^{-1}\cdot\Ad(U_{ij}^{++})\cdot(\rho_i\otimes
\rho_j^\o)\cdot\ga$.  By \cite{BE},1 Prop. 18, we have
$$[\be_{ij}]=
[\a^+_{\rho_i\otimes\id}]
[\a^-_{\id\otimes\rho_j^\o}]
= [\a^-_{\rho_j\otimes\id}]
[\a^+_{\rho_i\otimes\id}]
= [\a^-_{\id\otimes\rho_j^\o}]
[\a^+_{\rho_i\otimes\id}]$$
as $\M$-$\M$ sectors.

The following proposition is originally
due to Izumi \cite{I4} (with a different
proof) and first due to Ocneanu \cite{O6} in the setting
of the asymptotic inclusion.  (Also see \cite{EK3}.)

\begin{proposition}
Each $[\be_{ij}]$ is an irreducible $\M$-$\M$ sector and
these are mutually different for different pairs of
$i,j$.  Each irreducible
$\M$-$\M$ sector arising from $\N\otimes\N^\o\subset \M$ is of this
form.
\end{proposition}

\begin{proof}
We compute
$$\lan \be_{ij},\be_{kl} \ran=
\lan \a^+_{\rho_i\otimes\id}
\a^-_{\id\otimes\rho_j^\o},
\a^+_{\rho_k\otimes\id}
\a^-_{\id\otimes\rho_l^\o}\ran=
\lan \a^+_{\bar\rho_k\rho_i\otimes\id},
\a^-_{\id\otimes\rho_l^\o\bar\rho_j^\o}\ran.$$
The only sector which can be contained in
$[\a^+_{\bar\rho_k\rho_i\otimes\id}]$ and
$[\a^-_{\id\otimes\rho_l^\o\bar\rho_j^\o}]$ is the identity
by the above proposition.  So the above number is
$\de_{ik}\de_{jl}$.
Since the square sums of the statistical dimensions
for $\{\rho_i\otimes\rho_j^\o\}_{ij}$ and $\{\be_{ij}\}_{ij}$
are the same, it completes the proof.
\end{proof}

Note that here we have used the definition in \cite{LR} for the map
$\rho_i\otimes\rho_j^\o\mapsto \be_{ij}$, and a general theory
of this map has been studied in \cite{BE} under the name
$\a$-induction.   But in \cite{BE}, they assumed a certain
condition, called chiral locality in the terminology of \cite{BEK},
and some results in \cite{BE} depend on this assumption, while
the definition itself makes sense without it.  Our mixed use
of braidings $\e^+$ and $\e^-$ here violates this chiral locality
condition, so we can use the results in \cite{BE} here only when
they are independent of the chiral locality assumption.
For example, it is easy to see that the analogue of
\cite[I, Theorem 3.9]{BE} does not hold for our map here.

With the above proposition,
we have the following description of the dual
principal graph of $\N\otimes\N^\o\subset \M$ as a corollary,
which is originally due to Ocneanu \cite{O6}.
(Also see \cite{EK3}.)  Label even vertices with $(i,j)$ for
$[\be_{i\bar j}]$ and odd vertices with $k$
for $[\rho_k\otimes\id][\ga]$ and draw an edge with
multiplicity $N_{ij}^k$ between the even vertex $(i,j)$
and the odd vertex $k$.  The connected component of this
graph containing the vertex $(0,0)$ is the dual principal graph
of the subfactor $\N\otimes\N^\o\subset \M$, which is the same as
the principal graph.

We next study the tensor category of the $\M$-$\M$ sectors.

\begin{lemma}
Let $V, W$ be intertwiners from $\rho_i\rho_k$ to $\rho_m$
and from $\rho_j\rho_l$ to $\rho_n$, respectively,
in $\N$.
Then $V\otimes {W^*}^\o\in\N\otimes\N^\o$ in an intertwiner from
$\be_{ij}\be_{kl}$ to $\be_{mn}$.
\end{lemma}

\begin{proof}
By a direct computation.
\end{proof}

Then we easily get the following from the above lemma.
(The quantum $6j$-symbols for subfactors have been introduced
in \cite{O3} as a generalization for classical $6j$-symbols.
See \cite[Chapter 12]{EK2} for details.)

\begin{theorem}
In the above setting, the tensor categories of
$(\N\otimes\N^\o)$-$(\N\otimes\N^\o)$
sectors and $\M$-$\M$ sectors with quantum $6j$-symbols are
isomorphic.
\end{theorem}

\section{Relations with the quantum double}
\label{S5}
This section contains our main results.

Here below we will consider an inclusion $\A\subset \B$ of nets of 
factors.  We shall say that $\A\subset \B$ has finite index if there 
is a consistent family of conditional expectations $\E_I: \B(I)\to 
\A(I)$, $I\in\I$ and $[\B(I):\A(I)]_{\E_I}<\infty$ does not depend on 
$I\in\I$.  The independence of the index of the interval $I$ 
automatically holds if there is a vector (vacuum) with Reeh-Schlieder 
property and $\E_I$ preserves the vacuum state (standard nets, see 
\cite{LR}).  The index will be simply denoted by $[\B:\A]$.

\begin{proposition}\label{p15} Let $\A\subset \B$ be a finite-index
inclusion of nets of factors as above. If $\A$ and $\B$ are
completely rational then
$$
\mu_\A=\In^2\mu_\B
$$
with $\In=[\B:\A]$.
\end{proposition}

\proof If $\N_1,\N_2$ are factors, we shall use the symbol
$$
\N_1\stackrel{\alpha}{\perp} \N_2
$$
to indicate that $\N_1\subset \N'_2$ and $[\N'_2:\N_1]=\alpha$.

Let $E=I_1\cup I_2\in \I_2$; we will show that
$$
\begin{array}{ccc}\label{diag}
\B(E)&\stackrel{\mu_\B}{\perp}&\B(E')\\
\stackrel{\vphantom{\mu}}{\In^2}\,\cup&&\cup\ \In^2\\
\A(E)&\stackrel{\In^2\mu_\A}{\perp}&\A(E')
\end{array}
$$
where $\A(E)\subset \B(E)$ has index $\In^2$ because $\A(E)\isom
\A(I_1)\otimes \A(I_2)$, $\B(E)\isom \B(I_1)\otimes \B(I_2)$ and
$[\B(I_i):\A(I_i)]=\In$.

In the diagram, the commutants are taken in the Hilbert space $\H_\B$
of $\B$,
hence $\B(E)\stackrel{\mu_\B}{\perp}\B(E')$ is obvious.

We now show that on $\H_\B$
$$
\A(E)\stackrel {\In^2\mu_\A}{\perp} \A(E').
$$
Let $\gamma:\B\to \A$ be a canonical endomorphism with
$\lambda=\gamma\vert_\A$ localized in an interval $I_0$; then the net
$I\mapsto \A(I)$ on $\H_\B$ $(I\supset I_0)$ is unitarily equivalent to
the net
$$
I\mapsto \lambda(\A(I)) \quad\hbox{\rm on }\H_\A
$$
and we may assume $I_0\subset I_1$.

Then the correspondence associated with
$$
\A(E)\hbox{-}\A(E')\quad\hbox{\rm on }\H_\B\ ,
$$
namely $\H_\B$ with the natural commuting actions of $\A(E)$ and 
$\A(E')$, is unitarily equivalent to the one associated with
$$
\lambda(\A(E))\hbox{-}\lambda(\A(E'))\quad\hbox{\rm on }\H_\A\ ,
$$
namely $\H_\A$ with the commuting actions of $\A(E)$ and 
$\A(E')$ obtained by composing their defining actions with 
the map $X\to\lambda(X)$. But
$$
\lambda(\A(E))=\lambda(\A(I_1)\vee \A(I_2))=\lambda(\A(I_1))\vee
\A(I_2)
$$
and $\lambda(\A(E'))=\A(E')$ hence the $\A(E)$-$\A(E')$ correspondence
on $\H_\B$ is unitarily equivalent to
$$
(\lambda(\A(I_1))\vee \A(I_2))\hbox{-}\A(E')\quad\hbox{\rm on }\H_\A
$$
and its index is
$$
[\hat \A(E):\lambda(\A(I_1))\vee \A(I_2)]=[\hat
\A(E):\A(E)][\A(E):\lambda(\A(I_1))\vee \A(I_2)]=\mu_\A\In^2.
$$
It follows from the diagram that
$$
\In^2\mu_\A=\mu_\B\In^2\cdot \In^2,
$$
thus, $\In^2\mu_\B=\mu_\A$.
\endproof

The following Proposition may be generalized to the case of a
finite-index inclusion $\A\subset\B$ as above.

\begin{proposition}\label{c16bis} Let $\A$ be completely rational
with modular PCT and $\B\supset \A\otimes \A^\o$
be the LR net. Then also
$\B$ is completely rational with modular PCT.
\end{proposition}
\proof
Let $E=I_1\cup I_2$ and $I_3$ the bounded connected
component of $E'$. Set ${\cal C}\equiv \A\otimes \A^\o$.
Then the conditional
expectation  $\E_I:\B(I)\to{\cal C}(I)$ associated with the interval $I$,
where $I$ is the interior of
$\bar I_1\cup \bar I_2\cup \bar I_3$, maps $\B(E)$ onto
$\hat{\cal C}(E)$, because
$\E_I(\B(E))\subset {\cal C}(I_3)'\cap{\cal C}(I)=\hat{\cal C}(E)$, thus
\begin{equation}
        \E\equiv\E_0\cdot\E_I |_{\B(E)}
        \label{expect}
\end{equation}
is a finite-index expectation of $\B(E)$ onto ${\cal C}(E)$,
where $\E_0$ is the expectation of $\hat{\cal C}(E)$ onto ${\cal C}(E)$.
Therefore $\mu_{\B}<\infty$  follows by  a diagram similar to the one
in  (\ref{diag}) (with $\A\otimes\A^\o$ instead of $\A$),
as we know a priori that
the vertical inclusions have a finite index, while the bottom
horizontal inclusion has finite index by the argument given there.

Then the strong additivity of $\B$ follows easily, and so its modular 
PCT property, but we omit the arguments that are not essential here 
(if $\A$ is conformal case this follows directly because then also 
$\B$ is conformal).

We now show the split property of $\B$. For notational convenience
we treat the case of two separated intervals, rather than that of an
interval and the complement of a larger interval.
It will be enough
to show that the above expectation (\ref{expect}) satisfies
\[
\E(b_1 b_2) = \E(b_1)\E(b_2) \ ,\quad b_i\in{\cal C}(I_i) \ ,
\]
and $\E(\B(I_i))\subset {\cal C}(I_i)$, as we may then compose a normal
product state
$\phi_1\otimes\phi_2$ of ${\cal C}(I_1)\vee{\cal C}(I_2)\simeq
{\cal C}(I_1)\otimes{\cal C}(I_2)$ with $\E$ to get a normal product state of
${\cal B}(I_1)\vee{\cal B}(I_2)$.

Let $R_i^{(h)}\in\B(I_h)$, $h=1,2$, be elements satisfying
the relations (\ref{r}) so that $\B(I_h)$ is generated by ${\cal C}(I_h)$
and $\{R_i^{(h)}\}_i$. With $b_h\in\B(I_h)$ we then have
\[
b^{(h)}=\sum_i a_i^{(h)} R_i^{(h)} \ ,\quad a_i^{(h)} \in{\cal C}(I_h)\ ,
\]
hence
\[
b^{(1)}b^{(2)}=\sum_{i,j} a_i^{(1)} a_j^{(2)} R_i^{(1)} R_j^{(2)}\ ,
\]
so we have to show that $\E(R_i^{(1)} R_j^{(2)})=0$ unless
$i=j=0$. Now $R_i^{(1)}=u_i R_i^{(2)}$ for some unitary
$u_i\in\hat{\cal C}(E)$ and
\[
\E_I (R_i^{(2)} R_j^{(2)})=\E_I (\sum_k {C_{ij}^k}^{(2)}R_k^{(2)})
= {C_{ij}^0}^{(2)}= \delta_{\bar i j}{C_{ij}^0}^{(2)}\ ,
\]
(see Appendix \ref{A1} for the definition of the $C_{ij}^k$), hence
\begin{equation*}
\E (R_i^{(1)} R_j^{(2)})=\E (u_i R_i^{(2)} R_j^{(2)})
=\E_0 (u_i \E_I (R_i^{(2)} R_j^{(2)}))
= \E_0 (u_i     {C_{i\bar i}^0}^{(2)}) = \E_0 (u_i )
{C_{i\bar i}^0}^{(2)} \ ,
\end{equation*}
which is $0$ if $i\neq 0$ because $\E_0 (u_i )\in{\cal C}(E)$ is        
an intertwiner between  irreducible     endomorphisms localized
in $I_1$ and $I_2$, while $\E_0 (u_0) =\E_0 (1) = 1$.
\endproof

We get the following corollary, where the last part will
follow from Proposition \ref{p22} later.

\begin{corollary}\label{c16} Let $\A$ be completely rational  and
$$
\A\otimes \A^\o\subset \B
$$
be the LR inclusion. Then
$$
\mu^2_\A=\In^2_{\rm  global}\mu_\B
$$
where $\In_{\rm global}=\sum d(\rho_i)^2$.

In particular, $\mu_\B=1$ if and only if
$\A(E)\subset\hat \A(E)$ is isomorphic
to the LR inclusion.\end{corollary}
\proof By Propositions \ref{p15}, \ref{c16bis} and \ref{p22}.
\endproof
\begin{lemma}
\label{l17} Let $\A_1$, $\A_2$ be irreducible, Haag dual nets on 
separable Hilbert spaces.  Assume that each sector of $\A_1$ is of 
type I. If $\rho$ is an irreducible localized endomorphism of 
$\A_1\otimes \A_2$, then
$$
\rho\simeq\rho_1\otimes \rho_2
$$
with $\rho_i$ irreducible localized endomorphisms of $\A_i$.
\end{lemma}
\proof 
Let $\pi$ be a DHR representation of $\A_1\otimes \A_2$ (see 
Appendix \ref{rep-endo}) on a separable Hilbert space $\H$.  
Then $\pi(\mathfrak A_1)$ and $\pi(\mathfrak A_2)$ generate the
von Neumann algebra $\B(\H)$, where $\mathfrak A_i$ denotes the 
quasi-local C$^*$-algebra associated by $\A_i$.
Hence $\pi(\mathfrak A_1)''$ and $\pi(\mathfrak A_2)''$ are factors.  

Let $\pi_i\equiv\pi |_{\A_i}$, where we identify $\A_1$ with 
$\A_1\otimes\mathbb C$ and $\A_2$ with $\mathbb C\otimes\A_2$, then 
$\pi_i$ is easily seen to be localizable in bounded intervals (namely 
if $I_1\in\I$, the restriction of $\pi_1$ to the $C^*$-algebra 
generated by $\{ \A_i(I): I\in I'_1, I\in\I\}$ extends to a normal 
representation of $\A_i(I'_1)$).  Therefore $\pi_i$ is unitarily 
equivalent to a localized endomorphism of $\A_i$.  As $\pi_1$ is a 
factor representation, by assumption $\pi(\mathfrak A_1)''$ is a type 
I factor and so is $\pi(\mathfrak A_2)''$.  We then have a 
decomposition
$$
\pi=\pi_1\otimes \pi_2\ .
$$
This concludes the proof.
\endproof
\begin{corollary}\label{c18}
Let $\A$ be a completely rational net on a separable Hilbert space.  
The only irreducible finite dimensional sectors of $\A\otimes \A^\o$ 
are
$$
[\rho_i\otimes \rho_j^\o]
$$
with $[\rho_i]$, $[\rho_j]$ irreducible sectors of $\A$.
\end{corollary}
\proof Immediate by Lemma \ref{type I} and the above Lemma.
\endproof
\begin{lemma}\label{c19} Let $\A$ be  completely rational and $\B\supset
\A\otimes \A^\o$ the LR net.
If $\sigma$ is an irreducible localized endomorphism of $\B$
and $\sigma\prec\alpha^+_\rho$, $\sigma\prec\alpha^-_{\rho'}$
for some localized
endomorphism $\rho,\rho'$ of $\A\otimes \A^\o$, then $\sigma$ is
localized in a bounded interval.
\end{lemma}

\proof The thesis follows because $\sigma\prec\alpha^+_\rho$
is localized in a right half-line and $\sigma\prec\alpha^-_\rho$
in a left half-line.
\endproof

The following lemma extends Theorem \ref{intersect}.

\begin{lemma} \label{p18}
Let $\A$ be a completely rational  net, $\{[\rho_i]\}_i$ the
system of all irreducible sectors with finite dimension,
and $\B\supset \A\otimes \A^\o$ the LR net. The following are
equivalent:
\begin{description}
\item{$(i)$} The braiding of the net $\A$ is non-degenerate.
\item{$(ii)$} $\B$ has no non-trivial localized
endomorphism (localized in a bounded interval, finite index).
\end{description}
\end{lemma}

\proof  We use now an
argument in \cite{C}. Let $\sigma$ be a non-trivial
irreducible localized endomorphism of $\B$ localized in an interval,
with
$d(\sigma)<\infty$.

By Frobenius reciprocity
\begin{eqnarray*}
\sigma&\prec& \a^+_{\sigma^{\rm rest}},\\
\sigma&\prec& \a^-_{\sigma^{\rm rest}},
\end{eqnarray*}
where $\sigma^{\rm rest}=\gamma\cdot\sigma |_{\A\otimes \A^\o}$ and
$\gamma:\B\to
\A\otimes \A^\o$ is a canonical endomorphism.
Hence if $\rho_k\otimes\textnormal{id}\prec\sigma^{\rm rest}$ is an
irreducible sector with
$[\a^+_{\rho_k\otimes\textnormal{id}}]=[\sigma]$,
then by \cite{LR}, Prop. 3.9, the monodromy of
$\rho_k\otimes\textnormal{id}$
 with $\gamma |_{\A\otimes \A^\o}
=\sum\rho_i\otimes\rho_i^\o$ must be trivial, namely $\rho_k$ is a
non-trivial sector with degenerate braiding.

The converse is true, namely if $\rho_k$ is a non-trivial degenerate sector,
then $\a^+_{\rho_k\otimes\textnormal{id}}$ is a non-trivial sector of
$\B$ localized in a bounded interval.
\endproof

\begin{lemma}\label{t19'}
Let $\A$ be a completely rational net with modular PCT and let
$\{[\rho_i]\}_i$ be the system of all finite dimensional sectors of $\A$.
If $E=I_1\cup
I_2\in \I_2$, then
\[
\lambda_E = \bigoplus_i \rho_i \bar\rho_i|_{\A(E)},
\]
where $\lambda_E=\gamma_E |_{\A(E)}$, the $\rho_i$'s are localized in
$I_1$ and the $\bar\rho_i$'s are localized in
$I_2$
\end{lemma}
\proof
Let $j =\Ad J$, where $J$ is the modular
conjugation of $\A(0,\infty)$. Given $I\in \I$ we may identify
$\A(I)^\o$ with $j(\A(I))=\A(-I)$. We define a net $\tilde \A$ on $\R$
setting
\[
\tilde \A(I)\equiv \A(I)\otimes \A(I)^\o =  \A(I)\otimes \A(-I),
\quad I\in \I.
\]
With $I = (a,b)$ with $0<a<b$ and $E=I\cup -I$, let
$\gamma_E :\hat \A(E)\to \A(E)$ be the canonical endomorphism and
 $\lambda_E\equiv \gamma_E |_{\A(E)}$. We identify now $\lambda_E$
 with an endomorphism of $\eta_I$ of $\tilde \A(I)$ and want to show
that
 $\eta_I$ extends to a localized endomorphism of $\tilde \A$.

The proof is similar to the one of Theorem \ref{t13b}. With $d>c>b$,
by Lemma \ref{l14} there is an extension
$\eta$ of $\eta_{(a,b)}$ to $\tilde \A(a,d)$ with
$\eta |_{\tilde \A(b,d)}=\textnormal{id}$ and a canonical endomorphism
$\eta_{(a,d)}$
acting trivially on $\A(a,c)$
with a unitary $u\in\tilde \A(a,d)$ such that
\[
\eta = \Ad u \cdot \eta_{(a,d)}.
\]
Therefore $\Ad u |_{\tilde \A(-\infty,c)}$ is an extension of
$\eta_{(a,b)}$ to ${\tilde \A(-\infty,c)}$ which acts trivially on
$\tilde \A(-\infty,a)$ and on $\tilde \A(b,c)$. Letting
$c\to\infty$
we obtain the desired extension of $\eta_{(a,b)}$ to $\tilde \A$,
that
we still denote by $\eta$.

Now, by Lemma \ref{l17} for $\tilde \A$, every irreducible subsector of
$\eta$
will be equivalent to $\rho_h\otimes(j\cdot\rho_k\cdot j)$ for some
$h,k$, hence each irreducible subsector of $\lambda_E$ must be
equivalent to $\rho_h\cdot\bar\rho_k |_{\A(E)}$ where $\rho_h$ is
localized in $(a,b)$ and $\rho_k$ is
localized in $(-b,-a)$. By Theorem \ref{t11} this is possible if
and only if $h=k$.
\endproof

\begin{corollary}\label{t20}
Let $\A$ be completely rational with modular PCT. The following are
equivalent.
\begin{description}
\item{\rm (i)} The net $\A$ has no non-trivial sector with finite
dimension.
\item{\rm (ii)} The net $\A$ has no non-trivial  sector
(with finite or infinite dimension).
\item{\rm (iii)} $\mu_\A = 1$, namely  $\A(E)' = \A(E')$ for all
$E\in \I_2$.
\end{description}
\end{corollary}
\proof (i) $\Rightarrow$ (ii): It will be enough to show that every sector
(possibly with infinite dimension) $\rho$
of $\A$ contains the identity sector.
Given $E=I_1\cup I_2$ with $I_1, I_2\in\I$, we may suppose that $\rho$
is localized in $I_1$ and choose a sector $\rho'$ equivalent to $\rho$
and localized in $I_2$. If $u$ is a unitary with $\Ad u\cdot \rho
=\rho'$, then $u\in\hat \A(E)$, hence $u\in \A(E)$ by assumptions. Now
$\A(E)\simeq \A(I_1)\otimes \A(I_2)$ by the split property,
hence there exists a conditional expectation $\E : \A(E)\to
\A(I_{1})$ with $\E(u)\neq 0$, thus $\E(u)$ is a
non-zero intertwiner between $\rho$ and the identity.

(ii) $\Rightarrow$ (iii) follows by Lemma \ref{t19'}.

(iii) $\Rightarrow$ (i) follows by Th. \ref{t11} (or by Lemma
\ref{t19'}).
\endproof

The condition $\mu_\A =1$ is however compatible with the existence
of soliton sectors.

Note also that the condition that $\A(E)\subset \hat \A(E)$
has depth $\leq 2$ (equivalently
$\hat \A(E)$ is the crossed product of $\A(E)$ by a finite-dimensional
Hopf algebra) is equivalent to the innerness of the sector $\lambda$
extending $\lambda_E$
(because $\lambda_E$ is implemented by a Hilbert space of
isometries in $\hat \A(E)$ \cite{L3}), hence it is equivalent
to the the property that all irreducible sectors of $\A$ have
dimension 1 by Lemma \ref{t19'}.

The following  is the main result of this paper.

\begin{theorem}\label{c21} Let $\A$ be completely rational  with
modular PCT. Then
$$
\mu_\A=\In_{\rm global}\equiv
\sum d (\rho_i)^2
$$
and $\A(E)\subset\hat \A(E)$ is isomorphic to
the LR inclusion associated with $\A(I_1)\otimes
\A(I_2)$ and all the finite-dimensional irreducible sectors
$[\rho_i]$ of
$\A$.
\end{theorem}
\proof
$\hat \A(E)\supset \A(E)$ contains the LR inclusion by the following
Proposition \ref{p22}. Since $\mu_\A =\In_{\textnormal{global}}$ by
Lemma \ref{t19'} it has to coincide with the LR inclusion.
\endproof

\begin{corollary} Let $\A$ be completely rational and conformal.
The inclusions $\A(E)\subset\hat\A(E)$ are all isomorphic for $E\in\I_2$.
\end{corollary}
\proof 
If $I\in\I$ and the $\rho_i$'s are localized in $I$, for any 
given $I_1\in\I$ there is a M\"obius transformation giving rise to an 
isomorphism of $\A(I)$ with $\A(I_1)$ carrying the $\rho_i$'s to 
endomorphisms localized in $I_1$.  Therefore the isomorphism class of 
$\{\A(E),\lambda_E\}$ is independent of $E\in\I_2$.  Hence the LR 
inclusions based on that are isomorphic.  
\endproof 
Indeed, by using the uniqueness of the $III_1$ injective factor 
\cite{Co,H} and the classification of its finite depth subfactors 
\cite{P2} we have the following.
\begin{corollary}
Let $\A$ be completely rational and conformal.
The isomorphism class of the inclusion $\A(E)\subset\hat\A(E)$, 
$E\in\I_2$, dependes only on the tensor category of the sectors of 
$\A$, not on its model realization.
\end{corollary}
\proof
If $\A$ is non-trivial and $I$ is an interval, the $\A(I)$ is a 
$III_1$ factor and, as the split property hols, $\A(I)$ is injective 
(see e.g. \cite{L4}). Thus $\A(I)$ is the
unique injective $III_1$ factor \cite{H}.

By Popa's theorem \cite{P2}, if $\N$ is a $III_1$ injective factor and 
$\cal T\subset{\rm End}(\N)$ a rational tensor category isomorphic to 
the tensor category of sectors of $\A$ (as abstract tensor 
categories), then there exists an isomorphism of $\N$ with $\A(I)$ 
implementing the equivalence between the two tensor categories.

Since the LR inclusion $\N\otimes \N^{\o}\subset\M$ clearly depends, 
up to isomorphism, only on $\N$ and the tensor category $\cal 
T\subset{\rm End}(\N)$, it is then isomorphic to 
$\A(E)\subset\hat\A(E)$.
\endproof
We now show that, even in the infinite index case, the two-interval
inclusion always contains the LR inclusion associated with any
rational system of irreducible sectors.
\begin{proposition}\label{p22}
Let $\A$ be completely rational with modular PCT $j$ and $E=I\cup 
-I\in\I_2$ a symmetric 2-interval and $\{[\rho_i]\}$ a rational system 
of irreducible sectors of $\A$ with finite dimension, with the 
$\rho_i$'s localized in $I$.  Let $R_i\in (\hbox{\rm 
id},\bar\rho_i\rho_i)$ be non-zero intertwiners, where 
$\bar\rho_i=j\cdot\rho_i\cdot j$.

If $\M$ is the von Neumann subalgebra of $\hat A(E)$ generated by 
$\A(E)$ and $\{R_i\}_i$, then $\M\supset \A(E)$ is isomorphic to the 
LR inclusion associated with $\{[\rho_i]\}_i$, in particular
\[
[\M:\A(E)]=\sum_i d(\rho_i)^2.
\]
More generally this holds true if the assumption of complete rationality 
is relaxed with possibly $[\hat \A(E) : \A(E)]=\infty$.
\end{proposition}
\proof
Denoting by $\N$ the factor $\A(0,\infty)$, we may assume $\bar
I\subset (0,\infty)$ and consider the $\rho_i$ as endomorphisms of
$\N$.
Let then $V_i$ be the isometry standard implementation of $\rho_i$ as
in \cite{GL}.
Since $JV_i J = V_i$, we
have
$$
\rho_i\bar\rho_i(X)V_i=V_i X
$$
for all $X\in \N\vee \N'$, hence for all local operators $X$ by strong
additivity.

Since $\rho_i$ is irreducible, $({\rm id},\rho_i\bar\rho_i)$ is
one-dimensional, thus $R_i$ is a multiple of $V_i$ and we may assume
$R_i={\sqrt {d(\rho_i)}}V_i$,
thus
\begin{equation}
        R_i^*R_i =d(\rho_i).
        \label{r1}
\end{equation}
Now $V_iV_j$ is the standard implementation
of $\rho_i\rho_j$ on $\N$ hence by \cite[Proposition \A.4]{GL}, we
have
\begin{equation}
R_iR_j=\sum_k C^k_{ij}R_k,
        \label{r2}
\end{equation}
where
$C^k_{ij}$ is the canonical intertwiner between $\rho_k\bar\rho_k$ and
$\rho_i\rho_j\bar\rho_i\bar\rho_j$ given by
\begin{equation}
C^k_{ij}=\sum_h w_h j(w_h)\simeq \sum_h w_h\otimes j(w_h),
\label{r3}
\end{equation}
where the $w_h$'s form an orthonormal basis of isometries in
$(\rho_k,\rho_i\rho_j)$.

Setting $\rho_0 =\textnormal{id}$,  we also have
\begin{equation}
R_i^*=d(\rho_i)C^{0*}_{{\bar i}i} R_{\bar i}\label{r4}.
\end{equation}
Indeed the above equality holds up to sign by the $j$-invariance of
both members \cite[Lemma \A.3]{GL}, but the $-$ sign does not occur
because both members have positive expectation values on the vacuum
vector.

Now by the split property $\A(E)=\A(I)\vee \A(-I)\simeq \A(I)\otimes
\A(-I)$
and $\A(-I)=j(\A(I))$ can be identified with $\A(I)^\o$, therefore
$\M$ is isomorphic to the algebra generated by $\A(I)\otimes \A(I)^\o$
and multiple of isometries $R_i$ satisfying the above relations.
Moreover, there exists a conditional expectation from $\M$ to
$\A(I)\otimes \A(I)^\o$.

Corollary \ref{cor-char} then gives the desired isomorphism between 
$\A(E)\subset\M$ and the LR inclusion. (The Longo-Rehren inclusion
in \cite{Ma2}, as well as in \cite{LR}, is dual to the one in
this paper, but it does not matter here. Notice further that, in the 
conformal case, the $2$-interval inclusion $\A(E)\subset\hat\A(E)$ is 
manifestly self-dual.)

The above proof works also in the case case $\mu_\A=\infty$ thanks to
Prop. \ref{characterization}.
\endproof

\begin{corollary}\label{c23} Let $\A$ be completely rational with
modular PCT. Then the braiding of the tensor category of all sectors
of $\A$ is non-degenerate.
\end{corollary}

\proof With the notations in Corollary \ref{c16} we have
$\mu_\A^2=\In_{\textnormal{global}}^2\mu_\B$. On the other hand
$
\In_{\textnormal{global}}^2=\In_{\textnormal{global}}(\A\otimes
\A^\o)$,
hence
\[
\In_{\textnormal{global}}(\A\otimes
\A^\o)=\mu_\A^2=\In_{\textnormal{global}}^2\mu_\B
\]
therefore $\mu_\B =1$. By Corollary \ref{t20} we $\B$ has no
non-trivial
sector localized in a bounded interval and this is equivalent to the
non-degeneracy of the braiding by Lemma \ref{p18}.
\endproof

That $\mu_\A=\In_{\rm global}$ implies
the non-degeneracy of the braiding has been noticed in
\cite[Corollary 4.3]{Mu1}.

An immediate consequence of Corollary \ref{c23} follows from the
work \cite{R1}, where a model independent construction
of Verlinde's matrices $S$ and $T$ has been performed,
provided the braiding symmetry is
non-de\-ge\-ne\-rate, thus providing a corresponding representation
of the modular group $SL(2,\mathbb Z)$. Hence we have:

\begin{corollary} The Verlinde's matrices $T$ and $S$ constructed in
\cite{R1} are non-de\-ge\-ne\-ra\-te,
hence there exists an
associated representation of the modular group $SL(2,\mathbb Z)$.
\end{corollary}

\begin{corollary} Let $\A$ be completely rational with modular PCT.
Every sector of $\A$ is a direct sum of finite dimensional sectors.
\end{corollary}
\proof
Assuming the contrary, by Proposition \ref{inf} we have an irreducible
sector $[\rho]$ with infinite dimension. Let $E=I_1\cup I_2\in\I_2$
with $\rho$ localized in $I_1$ and $\rho'$ equivalent to $\rho$ and
localized in $I_2$. Let $u$ be a unitary in $(\rho,\rho')$. Then
$u\in\hat \A(E)$, hence it has a unique expansion
\[
u = \sum_i x_iR_i, \quad x_i \in \A(E),
\]
where $R_i$ are as in Proposition \ref{p22}. As
$xu=u\rho(x)$, $x\in \A(I_1)$, we have
\[
x\sum_i x_iR_i=
\sum_i x_iR_i\rho(x) = \sum_i x_i(\rho_i\cdot\bar\rho_i)(\rho(x))R_i
= \sum_i x_i\rho_i(\rho(x))R_i \quad\forall x\in \A(I_1),
\]
thus $x x_i= x_i\rho_i(\rho(x))$ for all $i$. As there is a $x_i\neq
0$,
by the split property  there is a non-zero intertwiner between
$\rho_i\cdot\rho$ and the identity. As $\rho_i$ and $\rho$ are
irreducible, this implies that $\rho$ is finite dimensional,
contradicting our assumption.
\endproof

\begin{corollary} Let $\A$ be conformal and completely rational. Then every representation 
on a separable Hilbert space is M\"obius covariant with positive energy.
\end{corollary}
\proof 
By the preceding result every such representation is a direct sum of irreducible sectors
with finite dimension. According to \cite{GL1} every finite dimensional sector is
covariant with positive energy, thus also a direct sum of such sectors.
\endproof

\section{$n$-Interval Inclusions}\label{S6}
\def\B{{\cal B}}  \def\H{{\cal H}}
\newcommand{\ol}{\overline}

In this section we extend the results on the 2-interval subfactors to
arbitrary multi-interval subfactors. Let $\A$ be a local, 
irreducible net on $S^1$. We assume $\A$  to be  completely rational
with modular PCT, so that our previous analysis applies.
Alternatively $\A$ may be assumed to be conformal
with $\mu_{\A}=[\hat\A(E):\A(E)]$ finite and independent of the
2-interval $E$; this setting will be needed to derive Cor. \ref{n-reg}.

If $E\in\I_n$ we set
$$
\mu_n=[\hat{\A}(E):\A(E)]\ .
$$
With this notation $\mu_\A=\mu_2$. We also consider the situation occurring
in representations different
from the vacuum representation: if
$\rho$ is a localizable representation
of  $\A$ (i. e. a DHR representation, that, on $S^1$, are just the
locally normal representations), we set
$\mu^\rho_n=[\rho(\A(E'))':\rho(\A(E))]$.
\begin{lemma} \label{m1}
$
\mu_n^\rho=\mu_1^\rho\,\mu_n \ , \quad \forall n\in\mathbb N.
$
\end{lemma}
\proof
Let $E=I_1\cup I_2\cup\dots\cup I_n\in\I_n$. We may suppose that
$\rho$ is an endomorphism of $\A$ localized in $I_1$.
Since $\rho$ acts trivially on $E'$, we have
$\rho(\A(E'))'=\A(E')'=\hat\A(E)$, thus the inclusion
$\rho(\A(E))\subset\rho(\A(E'))'$
is  a composition
\[
\rho(\A(E))\subset\A(E)\subset\rho(\A(E'))'=\hat\A(E)\ ;
\]
by the split property $\rho(\A(E))\subset\A(E)$ is isomorphic
to $\rho(\A(I_1))\otimes\A(I_2\cup\cdots\cup I_n)\subset
\A(I_1)\otimes\hat\A(I_2\cup\cdots\cup I_n)$, therefore
\[
\mu^{\rho}_n = [\hat\A(E):\A(E)]\cdot [\A(I_1):\rho(\A(I_1)]\ .
\]
\endproof
\begin{lemma} \label{sigma_n}
$\mu^{\rho}_n=d(\rho)^2 \, \mu_2^{n-1}
\ , \quad \forall n\in\mathbb N.$
\end{lemma}
\proof
By the index-statistics theorem \cite{L2}
we have $\mu_1^{\rho}=d(\rho)^2$, hence, by  Lemma \ref{m1},
we only need to show that $\mu_n=\mu_2^{n-1}$.
We proceed inductively. If $n=1$ the claim is trivially true.
Assume the claim for a given $n$ and let $E_{n}= I_1\cup\cdots\cup
I_n\in\I_n$ and $E_{n+1}= I_1\cup\cdots
\cup I_n\cup I_{n+1}\in\I_{n+1}$. Then
\[
\A(E_{n+1})=\A(E_n)\vee\A(I_{n+1})\subset
\hat\A(E_n)\vee\A(I_{n+1})\subset\hat\A(E_{n+1})\ ,
\]
thus, by the split property, $\mu_{n+1}=\mu_n \cdot
[\hat\A(E_{n+1}): \hat\A(E_n)\vee\A(I_{n+1})] \,$
and, by the inductive assumption, we have to show
that $\hat\A(E_n)\vee\A(I_{n+1})\subset\hat\A(E_{n+1})$ is equal to
$\mu_2$.
But the commutant of this latter inclusion
$\A(I'_{n+1})\cap\A(E'_n)\subset \A(E'_{n+1})$
has index is $\mu_2$ because, by the split property, turns out to be
isomorphic to $\A(I_{\ell}\cup I_r)\otimes\A(L)\supset
\hat\A(I_{\ell}\cup I_r)\otimes\A(L)$, namely
to a 2-interval inclusion tensored by a common factor, where
$I_{\ell}$ and $ I_r$ are the two intervals of $E'_{n+1}$ contiguous
to $I_{n+1}$ and $L$ is the remaining $(n-1)$-subinterval of $E'_{n+1}$.
\endproof
\begin{theorem}\label{mce} Let $\A$ be a local, irreducible
completely rational net with modular PCT.
Let $E=\cup_{i=1}^n I_i\in\I_n$ and
$\lambda^{(n)}=\gamma^{(n)}\vert \A(E)$ where $\gamma^{(n)}$ is a canonical
endomorphism from $\hat{\A}(E)$ into $\A(E)$. Then
\begin{equation}
 \lambda^{(n)}\cong\bigoplus_{i_1,\ldots,i_n} N_{i_1
 \ldots i_n}^0 \rho_{i_1}\rho_{i_2}\cdots\rho_{i_n}\ ,
   \label{rho_n}
\end{equation}
where $\{[\rho_i]\}_i$ are all the irreducible sectors with
finite statistics,
$\rho_{i_k}$ being localized in $I_k$.
$N_{i_1\ldots i_n}^0$ is the
multiplicity of the
identical endomorphism in the product $\rho_{i_1}\ldots\rho_{i_n}$.

The same results hold true if complete rationality is replaced by
conformal invariance and assuming $[\hat\A(E):\A(E)]=\In_{\rm global} <\infty$
independently of the 2-interval $E$.
\end{theorem}
\proof Let $I$ be an interval which contains $\cup_i I_i$ and let
$\rho_{i_k}$, $k=1,\ldots,n$, be irreducible
endomorphisms localized in $I_k$, respectively.
Then the intertwiner space between $\rho_{i_1}\rho_{i_2}\cdots
\rho_{i_n}$, considered as an
endomorphism of $\A(I)$, and the identity has dimension
$N_{i_1\ldots i_n}^0$. We are using here the equivalence between
local and global intertwiners, that holds either by strong additivity
of by conformal invariance \cite{GL}.
These intertwiners are multiples of isometries in $\hat{\A}(E)$.
Thus, by the argument leading to
Th. \ref{t11}, $\rho_{i_1}\rho_{i_2}\cdots\rho_{i_n}\vert_{\A(E)}$
is contained in $\lambda^{(n)}$ with
multiplicity $N_{i_1\ldots i_n}^0$. We have thus
proved the inclusion $\succ$ in (\ref{rho_n}).
Now the dimension of the endomorphism on the right hand side of
(\ref{rho_n}) has been
computed in \cite{X2}. For the sake of selfcontainedness we repeat the argument:
\begin{multline}
\sum_{i_1,\ldots,i_n} N_{i_1\ldots i_n}^0\,d({\rho_1})\cdots
   d({\rho_n}) =\sum_{i_1,\ldots,i_{n-1}} \left (\sum_{i_n}
   N_{i_1\ldots i_{n-1}}^{\ol{i_n}}\,d(\rho_{i_n})\right )
   d(\rho_{i_1})\cdots d(\rho_{i_{n-1}} )\\
     = \sum_{i_1,\ldots,i_{n-1}}\left (d({\rho_1})\cdots
   d(\rho_{i_{n-1}})\right )^2= \left (\sum_i d(\rho_i^2)\right )^{n-1} \ ,
\end{multline}
where we have used Frobenius reciprocity
$N_{i_1\ldots i_n}^0=
N_{i_1\ldots i_{n-1}}^{\ol{i_n}}$, the fact
$d(\rho)=d({\ol{\rho}})$ and the identity
$\sum_i \langle \rho_i,\rho\rangle d(\rho_i)=d(\rho)$.
On the other hand, we have
\[
d(\lambda^{(n)})=[\hat{\A}(E):\A(E)]=\mu_\A^{n-1}=\In_{\rm global}^{n-1}=
(\sum_i d(\rho_i)^2)^{n-1},
\]
where the first equality is obvious, the second is given by Lemma
\ref{sigma_n} and the last one follows from the results of the preceding
section. Thus the endomorphisms on
both sides of (\ref{rho_n}) have the same dimension, hence they are
equivalent.

The last claim in the statement follows by the same arguments and the
equivalence between local and global intertwiners.
\endproof
\begin{corollary} Let $\A$ be as in Th. \ref{mce}. If $E\in\I_n$, then
$\A(E)\subset\hat\A(E)$ is isomorphic to the $n$-th iterated LR
inclusion associated with $\N\equiv\A(I)$, $I\in\I$, and the system
of all sectors of $\A$ (considered as sectors of $\N$).

In particular, for a fixed $n\in\mathbb N$,
the isomorphism class of $\A(E)\subset\hat\A(E)$ depends only on the
superselection structure of $\A$ and not on $E\in\I_n$.
\end{corollary}
\proof
Let $E=I_1\cup\dots\cup I_n\in\I_n$ with $\bar E\subset
(0,\infty)$ and $n=2^k$. It follows by Lemma \ref{sigma_n}
and the split property that
\[
[\hat\A(E\cup -E):\hat\A(E)\vee \hat\A(-E)]= \In_{\rm global}\ .
\]
On the other hand, if the $\rho_i$'s are localized in $I_1$,
then the algebra generated by $\hat\A(E)\vee \hat\A(-E)$ and
the standard implementation isometries $V_i$ of
$\rho_i |_{\hat\A(E)}$ is the associated LR inclusion, analogously as
in Th. \ref{c21}, and is contained in $\hat\A(E\cup -E)$, hence
coincides with that by the equality of the indices.

The corollary then follows in the case $n=2^k$ by induction, once we
note that at each step the extension
 $\alpha^+_{\rho_i\otimes{\rm id}}$ from
 $\hat\A(E)\vee \hat\A(-E)$ to $\hat\A(E\cup -E)$
 is $\rho_i |_{\hat\A(E\cup -E)}$.

 The same is then true for an arbitrary $n$ by taking relative
 commutants.
\endproof

\section{Examples and further comments}\label{Examples}
Our results may be first illustrated by considering the case of an inclusion
of completely rational, local conformal irreducible
nets $\A\subset \B$, where
$\A=\B^G$ is the fixed-point of $\B$ with respect to the action of a
finite group $G$ and $\mu_\B=1$. Then $[\B:\A]=|G|$, thus by Prop.
\ref{p15},
$\In_{\textnormal{global}}(\A)=\mu_\A=|G|^2$. Now $\A$ has the DHR \cite{DHR}
irreducible sectors $[\rho_{\pi}]$ associated with $\pi\in\hat G$
and
\[
\sum_{\pi\in\hat G}d(\rho_{\pi})^2=|G|,
\]
therefore $\A$ has extra irreducible sectors $[\sigma_i]$
with
\[
\sum_i d(\sigma_i)^2= |G|^2 - |G|.
\]
For example, in the case of Ising model, we have
$\A= \B^{\Z_2}$ as above (but with $\B$ twisted local, yet this does
not alter our discussion), thus
$\mu_\A=4$ and thus $\sum d(\rho_i)^2 = 4$, so the standard
three sectors are the only irreducible sectors.
\smallskip

On the other hand, in the situation studied in \cite{Mu4},
the superselection category of $\A$ is equivalent to the
representation category of a twisted quantum double $D^\omega(G)$
with $\omega\in H^3(G,{\mathbb T})$. Since $D^\omega(G)$ is semisimple
we again have
\[
\sum_{\sigma\in\widehat{D^\omega(G)}}d(\sigma)^2=\dim\,D^\omega(G)=|G|^2=\mu
_\A.
\]
\smallskip

One may compare this with the situation occurring
on a higher dimensional spacetime.
There the strong additivity property may be replaced by the
requirement that $\A(\O'\cap\tilde\O)'\cap\A(\tilde \O)=\A(\O)$ if
$\O\subset\tilde\O$ are double cones. If $E\equiv \O_1\cup\O_2$, where
$\O_1$ and $\O_2$ are double cones with space-like separated closure,
the split property gives a natural isomorphism of $\A(\O_1)\vee\A(\O_2)$
with $\A(\O_1)\otimes\A(\O_2)$ and
\[
[\A(E')':\A(E)]=\In_{\rm global}=\sum_{\pi\in\hat
G}d(\rho_{\pi})^2 = |G|,
\]
where $G$ is the gauge group and the $\rho_{\pi}$'s are the DHR sectors
\cite{DHR} (there is no extra sectors). The reason for this difference
is that on $S^1$ the complement of a 2-interval is still a
2-interval, thus the inclusion $\A(E)\subset\hat\A(E)$ is self-dual,
while on the Minkowski spacetime the spacelike complement of
$\O_1\cup\O_2$
is a connected region producing no charge transfer inclusion.
\smallskip

The index $\mu_\A$ in the models given by
the loop group construction for $SU(n)_k$ has been computed in
\cite{X2}.
Our results apply in particular to these nets and the 2-interval
inclusion is  the LR inclusion associated with the corresponding
irreducible sectors
$\{[\rho_i]\}_i$.

We note  that in this case the 2-interval inclusion is {\it not}
the asymptotic inclusion of the corresponding
Jones-Wenzl subfactor \cite{J,We}, even up to tensoring
by a common injective III$_1$ factor.  Consider $SU(2)_k$
as an example.  The net has $k+1$ sectors and if we choose
the standard generator, we get a corresponding
subfactor of Jones with
principal graph $A_{k+1}$, up to tensoring a common
injective factor of type III$_1$, as in \cite{W}.
If we apply the construction of the asymptotic inclusion to
this subfactor, we get a ``quantum double'' of only the
sectors corresponding to the even vertices of $A_{k+1}$.
We get the same result, if we apply the LR construction
to the system of $\N$-$\N$ sectors (or $\M$-$\M$ sectors).
But the construction of a subfactor from 4 intervals gives
a ``quantum double'' of the system of
{\it all} the sectors, both even and odd.  If we want to
get this system from the asymptotic inclusion or the
Longo-Rehren inclusion, we have to use also bimodules/sectors
corresponding to the odd vertices of the (dual) principal graph.
In order to get this LR inclusion from
the construction of the asymptotic inclusion, we need to proceed
as follows.  Let $\{[\rho_i]\}_i$ be the set of all the
sectors for the net arising from the loop group construction
for $SU(n)_k$ as above.
Then for a fixed interval $I\subset S^1$, we consider
$(\bigoplus_i \rho_i)(\A(I))\subset \A(I)$ which has finite index
and finite depth.  Take a hyperfinite II$_1$ subfactor
$P\subset Q$ with the same higher relative commutants as
$(\bigoplus_i \rho_i)(\A(I))\subset \A(I)$.  Then the tensor
categories of the sectors with quantum $6j$-symbols of
$Q\vee(Q'\cap Q_\infty)\subset Q_\infty$ and $\A(E)\subset \hat \A(E)$
are isomorphic.
For this reason,
the index of the asymptotic inclusion of the Jones subfactor
with principal graph $A_{k+1}$
is half of that of the subfactor arising from 4 intervals
and the net for $SU(2)_k$.  For $SU(n)_k$, this ratio of the two
indices is $n$.
\smallskip

Finally we notice that there are models like the $SO(2N)_1$ WZW models,
see \cite{BO} or \cite{Mu4}, where all irreducible sectors have dimension
one, yet the superselection category $\2C$ is modular in agreement with our
results.
In these cases the fusion graph is disconnected, therefore the equivalent
categories of $\M-\M$ and of $\N\otimes\N^\o-\N\otimes\N^\o$ sectors are
proper subcategories of the categories $\2C\times\2C^\o\simeq\D(\2C)$, where
$\D(\2C)$ is the quantum double of $\2C$.
\smallskip

We close this section with a few questions.  Does there exist a net 
with only trivial sectors and non-trivial 2-interval inclusions (thus 
$\mu_\A = \infty$)?  Is strong additivity automatic in the definition 
of complete rationality?  Is the LR inclusion the only extension of 
$\N\otimes\N^\o$ with the given canonical endomorphism $\bigoplus_i 
\rho_i\otimes\rho_i^\o$?  
\appendix
\section{The crossed product structure of the LR inclusion}\label{A1}
Let $\N$ be an infinite factor and $\{[\rho_i]\}_i$ a rational system 
of irreducible sectors of $\N.$ The LR inclusion \cite{LR} is a 
canonical inclusion $\N\otimes\N^\o\subset\M$ associated with $\N$ and 
$\{[\rho_i]\}_i$ such that
\[
\lambda \simeq \bigoplus_i \rho_i\otimes\rho_i^\o \ ,
\]
where $\lambda$ is the restriction to $\N\otimes \N^\o$ of the
canonical endomorphism of $\M$ into $\N\otimes \N^\o$.

In \cite{LR} such an inclusion is obtained by a canonical choice of 
the intertwiners $T\in({\rm id},\lambda)$ and 
$S\in(\lambda,\lambda^2)$ that characterize the canonical endomorphism 
\cite{L3} (Q-system).  We now show the universality property of this 
inclusion and its crossed product structure, that will provide a 
different realization of it. By LR inclusion we will mean the upward 
LR inclusion.

We shall consider the free $^*$-algebra $\M_0$ generated by $\N\otimes 
\N^\o$ and elements $R_i$ satisfying the relations
\begin{equation}\label{r}
\left\{
\begin{array}{l} R_i x=(\rho_i\otimes\rho^{\o}_i)(x)R_i,
        \qquad x\in \N\otimes \N^\o \  , \\
     R_i^*R_i =d(\rho_i)\  ,\\
 R_iR_j=\sum_k C^k_{ij}R_k \  , \\
    R_i^*=d(\rho_i)C^{0*}_{{\bar i}i} R_{\bar i} \ ,
      \end{array}   \right.
\end{equation}
where $C^k_{ij}$ is the canonical intertwiner between 
$\rho_k\otimes\rho^\o_k$ and $\rho_i\rho_j\otimes\rho^\o_i\rho^\o_j$ 
given by $C^k_{ij}=\sum_h w_h\otimes j(w_h)$, with $j$ the antilinear 
isomorphism of $\N$ with $\N^\o$, and the $w_h$'s form an orthonormal 
basis of isometries in $(\rho_k,\rho_i\rho_j)$.

We equip $\M_0$ with the maximal C$^*$ semi-norm associated to
the representations of $\M_0$ whose restriction to $\N\otimes \N^{\o}$
are normal and denote by $\M$ the quotient of $\M_0$
modulo the ideal formed by the elements that are null with respect to
this seminorm and refer to $\M$ as the free reduced
pre-$C^*$-algebra generated by $\N\otimes \N^\o$ and the $R_i$'s.

\begin{proposition}\label{characterization} 
Let $\N$ be an infinite factor with separable predual and 
$\{[\rho_i]\}_i$ a rational system of finite-dimensional irreducible 
sectors of $\N$.

Let $\M$ be the free reduced pre-$C^*$-algebra generated by $\N\otimes 
\N^\o$ and elements $R_i$ satisfying the relations (\ref{r}) as above.

Then $\M$ is a factor and $\N\otimes \N^\o\subset \M$ is isomorphic to 
the LR inclusion associated with $\N$ and $\{[\rho_i]\}_i$.

In particular every element $X\in\M$ has a unique expansion
\[
X=\sum_i x_iR_i, \quad x_i\in\N\otimes\N^\o.
\]
\end{proposition}
In other words: if $\N\otimes \N^\o$ acts normally on a Hilbert space 
$\H$ and $R_i\in\B(\H)$ are elements satisfying the relations 
(\ref{r}), then the sub-algebra $\M$ of $B(\H)$ generated by 
$\N\otimes \N^\o$ and the $R_i$'s is a factor and $\N\otimes 
\N^\o\subset \M$ is isomorphic to the LR inclusion.
\proof
Clearly all elements of $\M$ have the form
\begin{equation}
        X=\sum_i x_iR_i, \quad x_i\in \N\otimes \N^\o,
        \label{module}
\end{equation}
and we may suppose that $\M$ acts on a Hilbert space so that $\N$ and
$\N^\o$ are weakly closed.

We now construct an conditional expectation $\E : \M\to \N\otimes 
\N^\o$.  Setting $\rho_0 = \textnormal{id}$, the expectation $\E$ may 
be defined by \begin{equation} \E(X)= x_0 \end{equation} for $X$ given 
by (\ref{module}), once we show that this is well-defined.  To this 
end we will apply the averaging argument in \cite{ILP}.

Let $\J$ be the set of all $x_0\in \N\otimes \N^{\o}$ such that
there exist $x_i\in \N\otimes \N^{\o}$, $i>0$, with $\sum_{i\geq 0}
x_i R_i
=0$. Clearly $\J$ is a two-sided ideal of $\N\otimes \N^{\o}$, hence
$\J=0$ (as we want to show) or $\J=\N\otimes \N^{\o}$
(we may suppose $\N$ to be of type III).
Suppose $\J\neq 0$ and let $X=1+\sum_{i> 0}x_i R_i=0$,
thus
$$
X=1+\sum_{i> 0}ux_i R_i u^* =
1+\sum_{i> 0}ux_i \rho_i\otimes\rho_i^\o(u^*)R_i= 0
$$
for all unitaries $u\in \N\otimes \N^{\o}$. Letting $u$ run in the
unitary group of a simple injective subfactor $\cal R$ of  $\N\otimes
\N^{\o}$
and taking a mean over this group, we have
$$
X=1+\sum_{i> 0}y_i R_i  = 0,
$$
where
$y_i\in \N\otimes \N^{\o}$ intertwines id and $\rho_i\otimes\rho_i^\o$
on $\cal R$, thus on all $\N\otimes \N^{\o}$ by the simplicity of
$\cal R$.
Since $\rho_i\otimes\rho_i^\o$ is irreducible, $y_i=0$, $i>0$,
and we have $1=0$, a contradiction.

Notice now that
\begin{equation*}
R_i R_i^*=d(\rho_i)R_i C^{0*}_{{\bar i} i}R_{\bar i} 
= d(\rho_i)\rho_i\otimes\rho_i^\o(C^{0*}_{{\bar i} i})R_i R_{\bar i}
=\sum_k d(\rho_i)\rho_i\otimes\rho_i^\o(C^{0*}_{{\bar i} i})
C^k_{i {\bar i}} R_k,
\end{equation*}
thus, by the conjugate equation in \cite{L2}, we have
\[
\E(R_i R_i^*)=d(\rho_i)\rho_i\otimes\rho_i^\o(C^{0*}_{{\bar i} i})
C^0_{i {\bar i}}=\frac{1}{d(\rho_i)}\ ,
\]
so every $X\in\M$ has the unique expansion
\begin{equation}\label{re}
X=\sum_i x_i R_i, \quad x_i=d(\rho_i)\E(XR_i^*)\ .
\end{equation}
Denoting by $\M_1\supset\N\otimes \N^\o$ the LR inclusion associated 
with $\N$ and $\{[\rho_i]\}_i$, $\M_1$ is generated by $\N\otimes 
\N^\o$ and elements $R_i'$, with an expectation $\E'$, satisfying the 
relations as in (\ref{r}) and (\ref{re}) \cite[Section 5]{Ma2}, hence 
the linear map 
\begin{equation}\label{iso}
\Phi: X\equiv\sum_i x_i R_i\in\M\to\Phi(X)\equiv\sum_i x_i R'_i \in\M_1 
\end{equation} 
is clearly a homorphism of $\M$ onto $\M_1$, which is the identity on 
$\N\otimes \N^\o$. $\Phi$ is clearly one-to-one by the uniqueness of 
the expansion (\ref{re}) both in $\M$ and in $\M_1$.
\endproof
Note that the above Proposition gives an alternative construction of the 
LR inclusion, which is similar to Popa's construction of the symmetric 
enveloping algebra \cite{P}, as follows.  Let $\N$ act standardly on 
$L^2(\N)$ and $V_i$ be the standard isometry implementing $\rho_i$.  
The $^*$-algebra $\mathfrak A$ generated by $\N$ and $\N'$ is 
naturally isomorphic to the algebraic tensor product 
$\N\circledcirc\N^\o$ and the operators $R_i\equiv 
\sqrt{d(\rho_i)}V_i$ satisfy the relations (\ref{r}) by \cite[Appendix 
A]{GL}. By the above argument there exists a 
conditional expectation $\E:\mathfrak B\to\mathfrak A$, where 
$\mathfrak B$ is the $^*$-algebra generated by $\mathfrak A$ and the 
$V_i$'s.  Taking a normal state $\phi$ of $\N$, the state 
$\tilde\phi\equiv\phi\circledcirc\phi^\o\cdot\E$ of $\mathfrak B$ 
gives by the GNS representation the LR inclusion 
$\pi_{\tilde\phi}(\mathfrak A)'' \subset\pi_{\tilde\phi}(\mathfrak 
B)''$ (Prop. \ref{characterization}).
\begin{corollary}\label{cor-char}
Let $\N$ be an infinite factor with separable predual and 
$\{[\rho_i]\}_i$ a rational system of finite-dimensional irreducible 
sectors of $\N$.

Let $\M$ be a von Neumann algebra with $\M\supset\N\otimes \N^{\o}$ 
and $R_i\in\M$ elements satisfying the relations (\ref{r}).  If $\M$ 
is generated by $\N\otimes \N^{\o}$ and the $R_i$'s, then $\N\otimes 
\N^{\o}\subset\M$ is isomorphic to the LR inclusion associated with 
$\{[\rho_i]\}_i$.

In particular $(\N\otimes \N^{\o})'\cap\M=\mathbb C$ and there exists 
a normal conditional expectation from $\M$ to $\N\otimes \N^{\o}$.
\end{corollary}
\proof
The proof is immediate, the isomorphism is obtained as in (\ref{iso}):
\[
X\in\M\to\sum_i d(\rho_i)\E(X R_i^*)R'_i \ ,
\]
(notations analogous to the ones in (\ref{iso}).
\endproof 
In the following we shall iterate the LR construction, in order to
describe the structure of multi-interval subfactors.

With $\N$ an infinite factor as above and $\{[\rho_i]\}_i$ a system
of irreducible sectors with unitary braiding symmetry, let $\alpha^+$ be the
induction map from sectors $\rho_i\otimes\rho^\o_j$ of
$\N\otimes\N^\o$ to sectors of the LR extension $\M_1\equiv\M$
defined by formula (\ref{induction}). Then
$\{\alpha^+_{\rho_i\otimes{\rm id}}\}_i$ is a system of irreducible
sectors of $\M$ with braiding symmetry and we may construct the
corresponding LR inclusion $\M_1\otimes\M_1^\o\subset \M_2$,
where the opposite of  $\alpha^+_{\rho_i\otimes{\rm id}}$
is $\alpha^+_{\bar\rho_i\otimes{\rm id}}$.
We may then iterate the procedure to obtain a tower $\M_1\subset
\M_2\subset\M_{2^k}\subset\cdots$ and thus an inclusion
\begin{equation*}
\N_n\subset\M_n\ ,\quad n=2^k\ ,        
\end{equation*}
where $\N_n\equiv
\N\otimes\N^\o\otimes\N\otimes\cdots\N\otimes\N^\o$
($2^k$ tensor factors). By construction
this inclusion has index $\In_{\rm
global}^{n-1}$ and we refer to it as the $n$-th iterated LR inclusion.
\begin{proposition}\label{iter}
Let $n=2^k$. The  $n$-th iterated LR inclusion
$\N_n\subset\M_n$ is irreducible. If $\gamma^{(n)}:\M_n\to\N_n$
is the canonical endomorphism, its restriction $\lambda^{(n)}=
\gamma^{(n)}|_{\N_n}$
is given by
\begin{equation}
\lambda^{(n)}\simeq\bigoplus_{i_1,i_2,\dots,i_n}
 N^0_{i_1 i_2\dots i_n}
\rho_{i_1}\otimes\rho^\o_{i_2}\otimes\cdots\otimes\rho^\o_{i_n}\ ,
        \label{ln}
\end{equation}
where $ N^0_{i_1 i_2\dots i_n}\equiv \langle {\rm id},
\rho_{i_1}\bar\rho_{i_2}\cdots\bar\rho_{i_n}\rangle$.
\end{proposition}
\proof By a computation similar to the one in Sect. \ref{S6},
$\lambda^{(n)}$ defined by formula (\ref{ln}) has dimension
\[
d(\lambda^{(n)})=\In^{n-1}_{\rm global}\ ,
\]
therefore the formula  $\lambda^{(n)}=\gamma^{(n)}|_{\N_n}$
will follow by showing that
$
\rho_{i_1}\otimes\rho^\o_{i_2}\otimes\cdots\otimes\rho^\o_{i_n}\prec
\gamma^{(n)}|_{\N_n}$ with multiplicity $N^0_{i_1 i_2\dots i_n}$
and this will also imply the irreducibility of
$\N_n\subset\M_n$ because then $\lambda^{(n)}\succ{\rm id}$
with multiplicity one.

But $\rho_{i_1}\otimes\rho^\o_{i_2}\otimes\cdots\otimes
\rho^\o_{i_n}$ is unitarily
equivalent to
$\rho_{i_1}\bar\rho_{i_2}\cdots\bar\rho_{i_n}\otimes{\rm
id}\otimes\cdots\otimes{\rm id}$ in $\M_n$, by applying iteratively
Lemma \ref{LL13}, hence we have the conclusion.
\endproof
Let now $m<n=2^k$ be an integer and set $\N_m$ be the alternate
tensor product of $k$ copies of $\N$ and $\N^\o$
\[
\N_m\equiv
\N\otimes\N^\o\otimes\N\otimes\cdots\N\otimes\N^\o\ ,
\quad m \ {\rm factors}.
\]
We then define the $m$-th iterated LR inclusion
\[
\N_m\subset \M_m\ ,
\]
where
$\M_m$ is defined as the relative commutant in $\M_n$ of the
remaining $n-m$ copies of $\N$ and $\N^\o$, i.e. $\M_{m}=
(\N'_m\cap\N_n)'\cap\M_n$. Note that
$\N_m\subset\M_m$ is an irreducible inclusion of factors because
$\N'_m\cap\M_m\subset\N'_n\cap\M_n=\mathbb C$.

Arguing similarly as above we then have:
\begin{proposition}
Proposition \ref{iter} holds true for all positive integer $n$
(in formula (\ref{ln}) $\rho^\o_{i_n}$ is $\rho_{i_n}$ if $n$ is odd).
\end{proposition}
\proof
Let $n=2^k$. Let
$\{ V_{i_1\dots i_n}^{\ell}: \ell =1,2,\dots  N_{i_1\dots
i_n}\}$ be a basis of isometries in the space of elements in $\M_n$
that intertwine
$\rho_{i_1}\otimes\rho^\o_{i_2}\cdots\otimes\rho^\o_{i_n}$
on $\N_n$.
Arguing as in Prop. \ref{characterization} we see that any element
$X\in\M_n$ has a unique expansion
\[
X=\sum_{i_1\dots i_n}
\sum_{\ell}x^{\ell}_{i_1\dots i_n}V^{\ell}_{i_1\dots i_n}\ ,\ \quad
x^{\ell}_{i_1\dots i_n}\in \M_n\ .
\]
Using this expansion it is easy to check that for $m<n$ the factor
$\M_m$ defined above is generate by $\N_m$ and the
$V^{\ell}_{i_1\dots i_n}$'s with $i_{m+1}=i_{m+2}=\cdots =i_n = 0$.
The rest then follows easily.
\endproof

\section{Nets on $\mathbb R$ and on $S^1$ and their representations.}
\label{rep-endo}
In our paper we deal with nets on $\mathbb R$, rather than nets on 
$S^1$, for various reasons: because this is the natural language for 
our arguments, because our results are valid for nets that are not 
necessarily conformal and, finally, because even if our analysis were 
restricted to conformal nets on $S^1$, our proofs would require 
the analysis more general nets on $\mathbb R$ (the $t=0$ LR net is not 
conformal).

In the next Section \ref{disint} we will however need to deal with 
nets on $S^1$ and their representations, and then conclude consequences
for nets on $\mathbb R$. Although the relations 
between nets on $\mathbb R$ and on $S^1$ and their representations is 
straightforward, we will describe explicitely this point here for the 
convenience of the reader. However, for simplicity, we consider only 
the case of strongly additive, Haag dual nets.

{\it Nets on $S^1$.} Let $\A$ be a net of von Neumann algebras on 
$S^1$ on a separable Hilbert space satisfying Haag duality.  We also 
assume the local von Neumann algebras $\A(I)$ to be properly infinite, 
which is automatically true if is the split property holds, or if $\A$ 
is conformal (except, of course, for the trivial net 
$\A(I)\equiv\mathbb C$).

A representation $\pi$ of $\A$ is, by definition, a map $I\in\I\to 
\pi_I$ that associates to each interval $I\in\I$ of $S^1$ a 
representation, on a fixed Hilbert space, of the von Neumann algebra 
$\A(I)$ such that $\pi_{\tilde I}|_{\A(I)}=\pi_I$ if $I\subset\tilde 
I$.  We shall say that $\pi$ is {\it locally normal} if $\pi_I$ is 
normal for all $I\in\I$ and that $\pi$ is {\it localizable} if $\pi_I$ 
is unitary equivalent to id$|_{\A(I)}$ for all $I\in\I$.  As the 
$\A(I)$'s are properly infinite the two notions coincide if $\pi$ acts 
on a separable Hilbert space.  Moreover every representation of $\A$ 
on a separable Hilbert space is automatically locally normal 
\cite{Tak}, thus localizable.

Denote by $C^*(\A)$ the universal $C^*$-algebra \cite{FRS} associated 
with $\A$ (see also \cite{GL1}).  For each $I\in\I$ there is a 
canonical embedding $\iota_I: \A(I)\to C^*(\A)$ and $\iota_{\tilde 
I}|_{\A(I)}=\iota_I$ if $I\subset\tilde I$; we identify $\A(I)$ with 
$\iota_I(\A(I))$ if no confusion arises.  There is a one-to-one 
correspondence between representations of the $C^*$-algebra $C^*(\A)$ 
and representations of the net $\A$, given by $\pi\to 
\{I\to\pi_I\equiv \p\cdot\iota_I\}$.  Locally normal representations 
of the net $\A$ correspond, of course, to locally normal 
representations of $C^*(\A)$.  We shall always assume our 
representations to act on a separable Hilbert space, thus local 
normality is automatic.

As Haag duality holds, a localizable representation $\pi$ of $C^*(\A)$ 
is unitarily equivalent to a representation of the form 
$\s_0\cdot\rho$, where $\s_0$ is the representation of $C^*(\A)$ 
corresponding of the identity representation of $\A$ (we shall however 
not need this result).

{\it Nets on $\mathbb R$.} Given a net $\A$ of von Neumann algebras on 
$S^1$ satisfying Haag duality we may associate a net $\A_0$ of Neumann 
algebras on $\mathbb R=S^1\setminus\{\infty\}$ (identification by 
Cayley transform) by setting
\[
\A_0(I)=\A(I)\ ,
\]
for all bounded intervals $I$ of $\mathbb R$. We call $\A_0$ the 
{\it restriction} of $\A$ to $\mathbb R$. Clearly, if $\A$ is 
strongly additive, then $\A_0$ is also strongly additive and 
satisfies Haag duality on $\mathbb R$ in the form
 \begin{equation}
 	\A(I)'=\A(\mathbb R\setminus I)\ ,
 	\label{Haagduality}
 \end{equation}
where $I\subset \mathbb R$ is either an interval or an half-line 
$(a,\infty)$ or $(-\infty, a)$, $a\in\mathbb R$.

Here, if $E\subset\mathbb R$ has 
non-empty interior, we denote by $\mathfrak A_0(E)$ the C$^*$-algebra 
generated by the von Neumann algebras $\A_0(I)$'s as $I$ runs in
the intervals contained in the region $E$ and set 
$\A_0(E)=\mathfrak A_0(E)''$.

Conversely, let now $\A_0$ be a strongly additive net of properly infinite 
von Neumann algebras $\A_0(I)$ on the (bounded, 
non-trivial) intervals of $\mathbb R$ satisfying Haag duality 
(\ref{Haagduality}). 

We may compactify $\mathbb R$ to $S^1=\mathbb
R\cup\{\infty\}$ and extend $\A_0$ to a net $\A$ on the intervals of
$S^1$ by defining
\begin{equation}\label{HD}
\A(I)\equiv\A_0(S^1\setminus I)'
\end{equation}
if $I$ is an interval whose closure contains the point $\infty$.  
Clearly, $\A$ is the unique Haag dual net on $S^1$ whose restriction 
to $\mathbb R$ is $\A_0$; we thus call $\A$ the {\it extension} of 
$\A_0$ to $S^1$.

We state explicitely this one-to-one in the following.
\begin{lemma}\label{corr}
Let $\A$ be a net on $S^1$ satisfying Haag duality
and strong additivity.  Then its restriction $\A_0$ to $\mathbb R$ 
satifies strong additivity and Haag duality on $\mathbb R$.

Conversely if $\A_0$ is a Haag dual (\ref{Haagduality}), strongly 
additive net on $\mathbb R$, then its extension $\A$ to $S^1$ is strongly 
additive and Haag dual.

Moreover $\A_0$ satisfies the split 
property on $\mathbb R$ if and only if $\A$ satisfies the split 
property on $S^1$.
\end{lemma}
\proof The proof is immediate. The statement concerning the split 
property follows because an inclusion of von Neumann algebras 
$\N\subset\M$ is split iff the commutant inclusion $\M'\subset\N'$ is 
split.
\endproof

We now consider the relation between representations of a net $\A$,
satisfying Haag duality and strong additivity on $S^1$ as in Lemma \ref{corr}
and its restriction $\A_0$ on $\mathbb R$.

A {\it DHR representation} $\pi_0$ of $\A_0$ is, by definition, a
representation $\pi_0$ of $\mathfrak A_0(\mathbb R)$ such that
$\pi_0|_{\mathfrak A_0(\mathbb R\setminus I)}$ is unitarily
equivalent to $\text{id}|_{\mathfrak A_0(\mathbb R\setminus I)}$
for every bounded non-trivial interval $I$ of $\mathbb R$, cf. 
\cite{DHR}.

Clearly a localizable representation $\pi$ of $\A$ determines a DHR 
representation $\pi_0$ of $\A_0$; indeed $\pi_0$ is consistently 
defined on $\cup_{a>0}\A(-a,a)$ by
\[
\pi_0(X)=\pi_I(X),\ X\in\A(I)\ ,
\]
where $I\equiv (-a,a)$, hence on all $\mathfrak A(\mathbb R)$ 
by continuity. We call $\pi_0$ the {\it restriction} of $\pi$ to $\A_0$.

Conversely, as we shall see,
every DHR representation  $\pi_0$ of $\mathfrak A_0(\mathbb R)$
determines uniquely a localizable representation $\pi$ of $\A$.

A {\it localized endomorphism} $\rho$ of $\A_0$ is, by definition, an 
endomorphism of $\mathfrak A_0(\mathbb R)$ such that $\rho|_{\mathfrak 
A_0(I')}=\text{id}|_{\mathfrak A_0(I')}$ for some interval 
$I\subset\mathbb R$; one then says that $\rho$ is localized in $I$.  
$\rho$ is transportable if for each interval $I_1$ there is an 
endomorphism $\rho_1$ localized in $I_1$ and (unitarily) equivalent to 
$\rho$ (as representations of $\mathfrak A_0(\mathbb R)$).  By Haag 
duality then $\rho_1 = {\rm Ad}u\cdot\rho$, where the unitary $u$ 
belongs to $\A_0(\tilde I)$, if $\tilde I$ is any interval containing 
both $I$ and $I_1$.  In this paper (as is often the case) {\it 
transportability is assumed} in the definition of localized 
endomorphism.

By a classical simple argument \cite{DHR}, a DHR representation 
$\pi_0$ of 
$\mathfrak A_0(\mathbb R)$ is unitarily
equivalent to a (transportable) endomorphism $\rho$ of 
$\mathfrak A_0(\mathbb R)$ localized in each given interval $I$; 
it is enough to put
\[
\rho(X)\equiv U\pi_0(X)U^*, \ X\in\mathfrak A_0(\mathbb R)\ ,
\]
where $U$ is a unitary intertwiner between 
$\pi_0|_{\mathfrak A_0(\mathbb R\setminus I)}$ and
${\rm id}|_{\mathfrak A_0(\mathbb R\setminus I)}$.

\begin{proposition}\label{corr2} Let $\A$ be a strongly additive, Haag 
dual net on $S^1$ and $\A_0$ be its 
restriction to $\mathbb R$, as in Lemma \ref{corr}.

If $\pi$ is a localizable representation of $\A$, its restriction 
$\pi_0$ to $\A_0$ is a DHR representation of $\A_0$.

Conversely, if $\pi_0$ is a DHR representation of $\A_0$, 
there exists a (obviously unique) localizable representation $\pi$ 
of $\A$ whose restriction to $\A_0$ is $\pi_0$.
\end{proposition}
\proof By the above discussion, we only show that if $\pi_0$
is a DHR representation of $\A_0$, 
there exists a localizable representation $\pi$ 
of $\A$ such that $\pi_I = \pi_0|_{\A(I)}$ if $I$ is a bounded 
interval of $\mathbb R$.

Indeed, if the closure of $I$ contains the point $\infty$, we can define
$\pi_I$ as the normal extension of $\pi_0|_{\mathfrak 
A_0(I\setminus\{\infty\})}$, once we show the necessary normality property.
Now the normality of 
$\pi_0|_{\mathfrak A_0(I\setminus\{\infty\})}$ does not depend on the 
unitary equivalence class of $\pi_0$, thus we may replace $\pi_0$ by a 
DHR endomorphism $\rho$ of $\A_0$ 
localized in interval $I_1\subset \mathbb R$ with $I_1\cap 
I=\emptyset$. But then $\rho|_{\mathfrak A_0(I\setminus\{\infty\})}$ 
is the identity, hence normal.
\endproof

By definition, the sectors of $\A$ (resp. of $\A_0$) are 
the unitary equivalence classes of localizable representations of $\A$ 
(resp. of DHR representations of $\A_0$). By the above discussions, the two classes 
are in one-to-one correspondence. 

On the other hand localizable representations of $\A$
corresponds to localizable representations of $C^*(\A)$ and
DHR representations of $\A_0$ are equivalent to DHR localized 
endomorphisms of $\A_0$, hence we have the following.

\begin{corollary}\label{corr3} Let $\A_0$ be a strongly additive, Haag 
dual as in (\ref{Haagduality}), net on $\mathbb R$ and $\A$ be its 
extension to $S^1$.  The restriction map $\pi\to\pi_0$ gives rise to a 
natural one-to-one correspondence between unitary equivalence classes 
of localizable representations of $C^*(\A)$ and unitary equivalence 
classes of DHR localized endomorphisms of $\A_0$.

In particular $\pi(C^*(\A))''=\pi_0(\mathfrak A_0(\mathbb R))''$, so 
$\pi$ is of type I iff $\pi_0$ is of type I.
\end{corollary}
\proof It remains to check the last part of the statement. As 
$C^*(\A)$ is generated (as C$^*$-algebra) by the von Neumann 
algebras $\A(I)$ as $I$ runs in the intervals of $S^1$, one has
$\pi(C^*(\A))''=\vee_I\pi_I(\A(I))$, thus clearly 
$\pi(C^*(\A))''\supset\pi_0(\mathfrak A_0(\mathbb R))''$. 

On the other hand if $I$ is an interval of $S^1$, by local normality 
and strong additivity we have 
$\pi_I(\A(I))=\pi_I(\A(I\setminus\{\infty\}))\subset\pi_0(\mathfrak 
A_0(\mathbb R))''$, hence $\pi(C^*(\A))''\subset\pi_0(\mathfrak 
A_0(\mathbb R))''$.  
\endproof
The naturality in the above corollary means that the tensor 
categories of localizable representations of $C^*(\A)$ and of DHR localized 
endomorphisms of $\A_0$ are equivalent, but we do not need this form 
of the above statement.

\section{Disintegration of locally normal representations
and of sectors.}\label{disint}
Takesaki and Winnink \cite{TW} have shown that a locally normal state
decomposes into locally normal states, if the split property holds. We
shall show here analogous results for localizable representations
(sectors). Our arguments work, however, along the same lines 
to show that locally normal representations decompose into locally
normal representations, also on higher dimensional manifolds.

We begin with a simple Lemma.
\begin{lemma}\label{KK} Let $\cal M$ be a von Neumann algebra,
$\mathfrak L\subset \cal M$ a $\s$-weakly dense C$^*$-subalgebra and
$J\subset\mathfrak L$ a right ideal of $\mathfrak L$.

If $\pi$ is  a representation of $\mathfrak L$ on a Hilbert space
$\H$ such that $\pi|_J$ is $\s$-weakly continuous and
$\overline{\pi(J)\H} =\H$, then
$\pi$ is $\s$-weakly continuous, thus it extends uniquely
to a normal representation of $\cal M$.
\end{lemma}
\proof
It is sufficient to show that $\pi$ is $\s$-weakly continuous on the
unit ball of $\mathfrak L$, see e.g. \cite{Tak}. Let then $\{a_i\}_i$ be a
bounded net of elements $a_i\in\mathfrak L$ such that $a_i\to 0$ $\s$-weakly.
If $t\in B(\H)$ is a $\s$-weak limit point of $\{\pi(a_i)\}_i$, we
have to show that $t=0$. By considering a subnet, if necessary, we may
assume $\pi(a_i)\to t$. Given $h\in J$, we have $a_i h\in J$ and
$a_i h\to 0$, thus $\pi(a_i h)\to 0$ because
$\pi|_J$ is $\s$-weakly continuous, therefore
\[
t\pi(h)= \lim_i \pi(a_i)\pi(h)= \lim_i \pi(a_i h) = 0,
\]
and this entails $t=0$ because $h$ is arbitrary and $\pi(J)\H$ is
dense in $\H$.
\endproof
We shall use the well-known fact that
the C$^*$-algebra of compact operators on a separable Hilbert space
$\H$ has only one non-degenerate (i.e. not
containing the zero representation) representation, up to multiplicity,
hence a unique normal extension to $B(\H)$.
\begin{corollary}\label{e} Let $\N$ be a type I factor with separable
predual, $K\subset\N$ the ideal of compact operator relative to
$\N$ and $\mathfrak L$  a C$^*$-algebra with
$K\subset\mathfrak L\subset \cal M$.

If $\pi$ is  a representation of $\mathfrak L$ such that $\pi|_K$ is
non-degenerate, then
$\pi$ is $\s$-weakly continuous, thus it extends uniquely
to a normal representation of $\N$.
\end{corollary}
\proof
Immediate because any non-degenerate representation of $K$ is
$\s$-weakly continuous and $K$ is $\s$-weakly dense in $\N$.
\endproof
Let $\A$ be a net of von Neumann algebras on $S^1$ over a 
separable Hilbert space satisfying the split property and Haag duality.

If $I,\tilde I$ are intervals, we write $I\subset\subset\tilde I$ if 
the closure of $I$ is contained in the interior of $\tilde I$.  For 
each pair of intervals $I\subset\subset\tilde I$ we choose an 
intermediate type I factor $\N(I,\tilde I)$ between $\A(I)$ and 
$\A(\tilde I)$ and let $K(I,\tilde I)$ be the compact operators of 
$\N(I,\tilde I)$ (there is a canonical choice for $\N(I,\tilde I)$ 
\cite{DL}, but this does not play a role here).  We denote by 
$\I_{\mathbb Q}$ the set of intervals with rational endpoints and by 
$\mathfrak A$ the $C^*$-subalgebra of $C^*(\A)$ generated by all 
$K(I,\tilde I)$ as $I\subset\subset \tilde I$ run in $\I_{\mathbb Q}$.  
Clearly $\mathfrak A$ is norm separable.

If $I_1\subset\subset\tilde I_1\subset I_2\subset\subset\tilde I_2$
then clearly $\N(I_1,\tilde I_1)\subset\N(I_2,\tilde I_2)$, but
$K(I_1,\tilde I_1)$ is not included in $K(I_2,\tilde I_2)$. For this
reason we define the C$^*$-algebras associated to pairs of intervals
$I\subset\subset\tilde I$
\[
\mathfrak L(I,\tilde I)\equiv \N(I,\tilde I)\cap \mathfrak A\ .
\]
As $\N(I,\tilde I)$ is the multiplier algebra of $K(I,\tilde I)$, 
$\mathfrak L(I,\tilde I)$ consists of elements of
$\mathfrak A$ that are multipliers of $K(I,\tilde I)$.

By definition
$K(I,\tilde I)\subset\mathfrak L(I,\tilde I)\subset \N(I,\tilde I)$
and $\mathfrak A$ is the $C^*$-subalgebra of $C^*(\A)$ generated by all
$\mathfrak L(I,\tilde I)$ as  $I\subset\subset \tilde I$ run in
$\I_{\mathbb Q}$.
\begin{lemma}\label{isot} If
$I_1\subset\subset\tilde I_1\subset I_2\subset\subset\tilde I_2$
are intervals then
\[
\mathfrak L(I_1,\tilde I_1)\subset\mathfrak L(I_2,\tilde I_2)\ .
\]
\end{lemma}
\proof
$\mathfrak L(I_1,\tilde I_1)\subset\N(I_1,\tilde I_1)
\subset\N(I_2,\tilde I_2)$, thus
\[
\mathfrak L(I_1,\tilde I_1)\subset\N(I_2,\tilde
I_2)\cap\mathfrak A = \mathfrak L(I_2,\tilde I_2)\ .
\]
\endproof
\begin{proposition}\label{non-deg}
Let $\pi$ be a locally normal representation of $C^*(\A)$. Then
$\pi|_{\mathfrak A}$ is a representation of $\mathfrak A$ and
$\pi|_{K(I,\tilde I)}$ is non-degenerate for every of pair of intervals
$I\subset\subset\tilde I$.

Conversely, if $\s$ is a representation of $\mathfrak A$ such that
$\s|_{K(I,\tilde I)}$ is non-degenerate for
all intervals $I,\tilde I\in \I_{\mathbb Q}$, $I\subset\subset\tilde
I$, there exists a unique locally normal
representation $\tilde\s$ of $C^*(\A)$ that extends $\s$.

Moreover equivalent representations $C^*(\A)$ correspond to equivalent
representations of $\mathfrak A$.
\end{proposition}
\proof
The only non-trivial part is that $\s$ extends to a locally
normal representation $\tilde\s$ of $C^*(\A)$.
If $I\subset\subset\tilde I$ are intervals in $\I_{\mathbb Q}$,
we denote by $\tilde\s_{I,\tilde I}$ the unique normal extension of
$\s|_{\mathfrak L(I,\tilde I)}$ to $\N(I,\tilde I)$ given by Corollary
\ref{e}.

Given an interval $I$,
we choose $I_1,\tilde I_1\in\I_{\mathbb Q}$,
$I_1\subset\subset \tilde I_1$  such that $I\subset\subset I_1$ and set
\[
\tilde\s_I\equiv\tilde\s_{I_1,\tilde I_1}|_{\A(I)}\ ,
\]
We have to show that $\tilde\s_I$ is well-defined, then $I\to \tilde\s_I$
is clearly a representation of $\A$.

Indeed, let $I_2,\tilde I_2\in\I_{\mathbb Q}$ with
$I_2\subset\subset\tilde I_2$ be another pair such that
$I\subset\subset I_2$.
We can choose $I_3,\tilde I_3\in\I_{\mathbb Q}$ such that
$I\subset\subset I_3\subset\subset\tilde I_3\subset\subset I_1\cap I_2$.
Then by Lemma \ref{isot}
$\mathfrak L(I_3,\tilde I_3)\subset\mathfrak L(I_i,\tilde I_i)$, $i=1,2$,
and therefore
\[
\tilde\s_{I_3,\tilde I_3}=\tilde\s_{I_1,\tilde I_1}|_{\N(I_3,\tilde I_3)}
=\tilde\s_{I_2,\tilde I_2}|_{\N(I_3,\tilde I_3)}\ .
\]
This concludes the proof.
\endproof
\begin{proposition}\label{ex}
Let $\pi$ be a locally normal representation of $C^*(\A)$ on a
separable Hilbert space and denote by $\pi_{\mathfrak A}$
be the restriction of $\pi$ to $\mathfrak A$. If
\[
\pi_{\mathfrak A}=\int_{X}^{\oplus}\!\! \pi_\lambda d\mu(\lambda)
\]
is a decomposition into irreducible representations $\pi_{\lambda}$
(which always exists), then $\pi_\lambda$ extends to a locally normal
representation $\tilde\pi_\lambda$ of $C^*(\A)$ for almost all
$\lambda$.
\end{proposition}
\proof By Proposition \ref{non-deg}, it is sufficient to show that
there exists a null set $E\subset X$ such that
$\pi_\lambda|_{K(I,\tilde I)}$ is non-degenerate for $\lambda\notin E$
and all $I,\tilde I\in\I_{\mathbb Q}$ with
$I\subset\subset\tilde I$. This is clear for a fixed pair $I,\tilde
I$ of the family, because $\pi_{K(I,\tilde I)}$ is non-degenerate. Then the
statement follows since the considered family of $K(I,\tilde I)$'s is countable.
\endproof
\begin{proposition} With the notations in the Proposition \ref{ex},
if $\pi(C^*(\A))''$ is a factor
not of type I, then for each $\lambda\in X$ the set
$X_{\lambda}\equiv
\{\lambda'\in X,
\pi_{\lambda'}\simeq \pi_{\lambda}\}$ has measure zero.
\end{proposition}
\proof The set $X_{\lambda}$ is measurable by Lemma \ref{measurable}
below.
We have $\mu(X\setminus X_{\lambda})>0$, as otherwise
$\pi$ would be quasi-equivalent to $\pi_{\lambda}$, hence
$\pi(\mathfrak A)''$
would be a type I factor. If $\mu(X_{\lambda})>0$, then
$\pi_{\mathfrak A}$ would be
the direct sum of two inequivalent representations
\[
\pi_{\mathfrak A}=\int_{X_{\lambda}}^{\oplus}\!\!\! \pi_\lambda
d\mu(\lambda) \oplus \int_{X\setminus X_{\lambda}}^{\oplus}\!\!
\pi_\lambda
d\mu(\lambda)
\]
which is not possible since $\pi(\mathfrak A)''$ is a factor.
\endproof

\begin{corollary}\label{inf2} If there exists a localizable
representation $\pi$ of $C^*(\A)$
with $\pi(C^*(\A))''$ a factor not of type I, then there exist
uncountably many inequivalent irreducible localizable representations of
$C^*(\A)$.
\end{corollary}
\proof
If the representation $\pi$ is factorial not of type I, then
the family of the $\pi_{\lambda}$'s
in the above proposition contains an uncountable set of
mutually  inequivalent irreducible localizable representations as desired.
\endproof

\begin{corollary}\label{inf} Let $\A_0$ be a strongly additive, split 
net of von Neumann algebras on the intervals of $\mathbb R$ which is Haag 
dual as in (\ref{Haagduality}).
If there exists a DHR localized endomorphism $\rho$ of $\A_0$
with $\rho(\mathfrak A_0(\mathbb R))''$ a factor not of type I, then there exist
uncountably many inequivalent irreducible DHR localized endomorphisms 
of $\A_0$.
\end{corollary}
\proof
Immediate by Corollary \ref{inf2} and Corollary \ref{corr3}.
\endproof

Before concluding this appendix we have to prove a Lemma 
that has been used. 
Let $\mathfrak A$ be any separable $C^*$-algebra and $\sigma$ a
representation of $\mathfrak A$. Choose a sequence of elements
$a_{\ell}\in\mathfrak A$ dense
in the unit ball $\mathfrak A_1$, a sequence $\phi_i\in \mathfrak
A^*$ dense in the Banach space of
normal linear functionals $(\sigma(\mathfrak A)'')_*$ associated with
$\sigma$.
A linear functional
$\phi\in \mathfrak A^*$ is then normal with respect to $\sigma$ if
and only if
\begin{equation}
        \forall k\in\mathbb N,\ \exists i\in\mathbb N :
        |\phi(a_{\ell})-\phi_i(a_{\ell})|\leq \frac{1}{k},\ \forall
        \ell\in\mathbb N.
        \label{approx}
\end{equation}
We thus have the following.

\begin{lemma}\label{measurable} Let $\mathfrak A$ be a separable
$C^*$-algebra, $\pi$ a representation of $\mathfrak A$ on a separable
Hilbert space and $\pi =
\int_X^{\oplus}\pi_{\lambda}d\mu(\lambda)$ a direct integral
decomposition into a.e.
irreducible representations $\pi_{\lambda}$ of $\mathfrak A$.
For any irreducible representation
$\sigma$ of $\mathfrak A$, the set $X_\sigma\equiv\{\lambda,
\pi_{\lambda}\simeq\sigma\}$ is measurable.
\end{lemma}
\proof
Let $\xi=\int_X^{\oplus}\xi(\lambda)d\mu(\lambda)$ be a vector with
$\xi(\lambda)\neq 0$, for all $\lambda\in X$, and consider the
functional of $\mathfrak A$ given by
$\phi_{\lambda}=(\pi_{\lambda}(\cdot)\xi(\lambda),\xi(\lambda))$.

As both $\sigma$ and $\pi_{\lambda}$ are irreducible, we have
$\sigma\simeq\pi_{\lambda}$ if and only if $\phi_{\lambda}$ is normal with
respect
to $\sigma$. With the previous notations, we then have by eq.
(\ref{approx})
\[
X_{\sigma}=\bigcap_{k}\bigcup_i\bigcap_{\ell}X_{ik\ell}
\]
where
\[
X_{ik\ell}=\{\lambda\in X:
        |\phi_{\lambda}(a_{\ell})-\phi_i(a_{\ell})|\leq \frac{1}{k}\}.
        \]
As      $X_{ik\ell}$ is measurable, also $X_{\sigma}$ is measurable.
\endproof
\medskip

\noindent{\bf Acknowledgments.}
A part of this work was done during visits of the first-named
author to Universit\`a di Roma ``Tor Vergata''.  Y.K. acknowledges
the hospitality and financial supports of CNR (Italy), Universit\`a
di Roma ``Tor Vergata'' and the Kanagawa Academy of Science and
Technology Research Grants.
R.L. wishes to thank the Japan Society for the Promotion
of Science for the invitation at the University of Tokyo in June
1997. The authors would like to thank K.-H.~Rehren for comments.

\end{document}